\documentclass[10pt]{article}

\usepackage{amssymb,amsmath,amsthm,verbatim}
\usepackage{graphicx}
\usepackage{lscape}
\usepackage{indentfirst}
\usepackage{latexsym}
\usepackage{multirow}
\usepackage{multicol}
\usepackage{longtable}
\usepackage{amssymb,amsmath,amsthm,verbatim, graphics}
\newtheorem{theorem}{Theorem}[section]

\newtheorem{lemma}[theorem]{Lemma}
\newtheorem{prop}[theorem]{Proposition}
\newtheorem{cor}[theorem]{Corollary}

\newtheorem{defn}{Definition}[section]

\newcommand{\bs}{\backslash}
\newcommand{\Z}{\mathbb{Z}}
\newcommand{\Q}{\mathbb{Q}}
\newcommand{\ds}{\frac{\partial}{\partial s}}
\newcommand{\Af}{\mathbb{A}_f}
\newcommand{\A}{\mathbb{A}}
\newcommand{\C}{\mathbb{C}}
\newcommand{\h}{\mathfrak{h}}
\newcommand{\R}{\mathbb{R}}
\newcommand{\F}{\mathbb{F}}
\newcommand{\K}{{\cal K}}
\newcommand{\vpi}{v}
\newcommand{\RK}{{\cal R}}
\newcommand{\D}{\mathfrak{D}}
\newcommand{\Pro}{\mathbb{P}}
\newcommand{\Po}{\mathcal{P}}
\newcommand{\SL}{\text{SL}_2}
\newcommand{\GL}{\text{GL}_2}
\newcommand{\M}{\text{M}_2}
\newcommand{\tr}{\text{tr}}
\newcommand{\Oh}{\mathcal{O}}
\newcommand{\X}{\mathcal{X}}
\newcommand{\Xs}{\mathcal{X}^*}
\newcommand{\no}{\text{n}}
\newcommand{\No}{\text{N}}
\newcommand{\Div}{\text{div}}
\font\cute=cmitt10 at 12pt

\newcommand{\kay}{{\text{\cute k}}}

\newcommand{\E}{\mathcal{E}}
\newcommand{\q}{{\bf q}}
\newcommand{\Ord}{\text{Ord}}
\newcommand{\ord}{\text{ord}}

\newcommand{\Hom}{\text{Hom}}
\newcommand{\cond}{\text{cond}}

\newcommand{\MSL}{\widetilde{\text{SL}_2}}
\newcommand{\sign}{\text{sign}}
\newcommand{\expo}{{\bf e}}
\newcommand{\Weil}{\rho_{\Lambda_L}}
\newcommand{\Weilb}{\overline{\rho}_{\Lambda_L}}
\newcommand{\Go}{\Gamma_0}
\newcommand{\MGo}{\widetilde{\Gamma_0}}

\begin{document}
\title{Singular Moduli of Shimura Curves} \author{Eric Errthum} \maketitle
\begin{abstract}
The $j$-function acts as a parametrization of the classical  
modular curve. Its values at complex multiplication (CM) points are called singular moduli and are algebraic integers. A 
Shimura curve is a generalization of the modular curve and, if the Shimura 
curve has genus $0$, a rational parameterizing function exists and when evaluated at a CM point 
is again algebraic over $\Q$. This paper shows that the coordinate maps given in \cite{Elk} for the Shimura curves associated to the quaternion algebras with discriminants $6$ and $10$ are Borcherds lifts of vector-valued 
modular forms. This property is then used to explicitly compute the rational norms of singular moduli on these curves.
This method not only verifies the conjectural values for the rational CM points listed in \cite{Elk}, but also provides a way of algebraically calculating the norms of CM points with arbitrarily large negative discriminant.
\end{abstract}

\begin{section}{Introduction}

The classical modular curve $\Xs_1$ is given as  the one-point compactification of the Riemann surface $\GL(\Z)\bs\h^\pm$ where $h^\pm = \Pro^1(\C)-\Pro^1(\R)$. Since $\Xs_1$ is a genus-$0$ surface, there exists an isomorphism $\Xs_1 \stackrel{\sim}{\rightarrow} \Pro^1$. The classical choice of such a map has Fourier expansion 
\begin{eqnarray*}
j(\tau) = \frac{1}{\q}+744+196884\q + \dots \in \frac{1}{\q}\Z[[\q]],
\end{eqnarray*}
(where $\q = e^{2\pi i \tau}$) at the cusp at $\infty$. The $j$-function also provides an identification of points on the modular curve with isomorphism classes of elliptic curves. When the associated elliptic curve has an extra endomorphism called complex multiplication (CM), $\tau$ is an irrational quadratic imaginary point of $\h^\pm$ and is called a CM point. A singular modulus is a value of the $j$-function at a CM point and is an algebraic integer. In 1984, Gross and Zagier \cite{GZ} gave an explicit formula to compute the norms of singular moduli. 

A Shimura curve is a generalization of the modular curve. Let $B$ be the quaternion algebra over $\Q$ with discriminant $D=D(B) > 1$ and let $\Gamma^* = N_{B^\times}(\Oh) \subset B^\times$ be the normalizer of a maximal order $\Oh \subset B$. Since there is an algebra embedding $B \hookrightarrow \M(\R)$, the discrete group $\Gamma^*$ embeds into $\GL(\R)$ and hence acts on $\h^\pm$. The Shimura curve $\Xs_D$ is then given as 
$$\Xs_D = \Gamma^* \bs \h^\pm.$$ 
When $B$ is a division algebra, $\Xs_D$ is a compact Riemann surface without cusps. 

Points on a Shimura curve can also be identified with certain 2-dimesional abelian varieties and again there is the notion of CM points. 
As before, there will be a generator of the function field, or Hauptmodul, $t_D: \Xs_D \rightarrow \Pro^1$, and, if properly normalized, the image of a CM point under $t_D$ will be algebraic over $\Q$. However, since $\Xs_D$ has no cusps, such a map does not have a $\q$-expansion and example calculations are more difficult than in the classical case. 

In \cite{Elk}, Elkies considered the cases of $D=6$ and $D=10$. First, by identifying which quadratic imaginary fields have class group $(\Z/2\Z)^r$ for $r\le2$, he determined which CM points have rational coordinates on $\Xs_D$.  Then 
with $\Gamma^*(l) = \{\gamma \in \Gamma^* \mid \gamma \equiv 1 \mod l\}$,
Elkies used explicit calculations of the geometric involution on $\Xs_D(l)= \Gamma^*(l) \bs \h^\pm$ for small primes $l$ to compute the coordinates for about half of the rational CM points on $\Xs_6$ and $\Xs_{10}$. 
The involutions on $\Xs_D(l)$ for higher $l$ are unknown and are needed to explicitly find the coordinates of the remaining half of the CM points using this method. Elkies does, however, provide a table of conjectural values for the remaining CM points obtained via numerical approximations and their behavior under standard transformations. 

In this paper, we use an alternate method that arises out of the theory of Borcherds forms to calculate the norms of singular moduli on the Shimura curves $\Xs_6$ and $\Xs_{10}$ and, as a special case, algebraically prove the conjectural values listed in \cite{Elk}. 
Although the methods are only demonstrated here for $D=6$ and $D=10$, the techniques should extend to a larger class of functions $\Xs_D \rightarrow \Pro^1$ for arbitrary indefinite discriminants $D$.

Let $L$ be a lattice in a rational inner product space $V \subset B$ with signature $(n,2)$ and let $L^\vee$ be its integral dual. Then a meromorphic modular form $F$ valued in $\C[L^\vee/L]$ can be given by its Fourier expansion 
\begin{eqnarray}
F(\tau) = \sum_{\eta \in L^\vee\!/\!L}\ \sum_{m\in \Q}c_\eta(m)\q^me_\eta, \label{FdefIntro}
\end{eqnarray}
where $e_\eta$ is the basis element of $\C[L^\vee/L]$ corresponding to $\eta$. When $c_\eta(m) \in \Z$ for $m<0$, $c_0(0)=0$, and $F$ has weight $1-\frac{n}{2}$, Borcherds \cite{BorGrass} constructs a form $\Psi(F): \Xs_D \rightarrow \Pro^1$ and gives its divisor in terms of rational quadratic divisors weighted by the coefficients $c_\eta(m)$ for $m<0$. In this more general setting $\Xs_D$ is formed by $B^\times$ acting on the product of the adeles of $B$ (viewed as an algebraic group) modulo a compact open set and a space of oriented negative 2-planes arising from the inner product.

Recently, Schofer \cite{JS} provided an explicit formula in terms of the coefficients of Eisenstein series for the norm
\begin{eqnarray}
\prod_{z \in Z(\Delta)}||\Psi(z,F)||^2 \label{CMProdIntro}
\end{eqnarray}
where $Z(\Delta)$ is the set of CM points of discriminant $\Delta$ on $\Xs_D$. As a corollary, he showed that since the $j$-function was in fact a Borcherds form, the Gross-Zagier factorization of singular moduli was a specific case of his main theorem.

In the cases of $D=6$ and $D=10$, the coordinate map $t_D : \Xs_D \stackrel{\sim}{\rightarrow} \Pro^1$ given in \cite{Elk} is defined by its divisor and normalized by its value at a chosen point. We show how this divisor can be expressed in terms of rational quadratic divisors. We then find a meromorphic modular form $F_D$ as in (\ref{FdefIntro}) that satisfies $\Div(\Psi(F_D)^2) = \Div(t_D)$. In the cases analyzed here, $n=1$ and the lattice $L$ arises as the trace-zero elements of $\Oh$. Then the proper vector-valued form $F_D$ is lifted from a scalar-valued modular form that is a linear combination of Dedekind-$\eta$ products. Next we compute a normalization constant, $c_D$, by applying (\ref{CMProdIntro}) to a base case. Since the divisors are equal and the two functions agree on the chosen base point, we conclude
$$\Psi(F_D)^2 = c_D t_D.$$

Finally, (\ref{CMProdIntro}) is used to calculate the norm of any CM point on $\Xs_D$. Since this method is a general calculation of norms, the tables of rational CM points found in \cite{Elk} arise as specific cases. For example, we can recompute known values, e.g. $$t_6(\Po_{-147}) = -\frac{11^423^4}{2^{10}3^35^67},$$  but can also explicitly verify the conjectural values such as $$ t_6(\Po_{-163}) = \frac{3^{11}7^419^423^4}{2^{10}5^611^617^6}.$$ Moreover, we can algebraically compute the norm of CM points with arbitrarily large discriminants. For example $t_6(\Po_{-996})$ is an algebraic number of degree 6 over $\Q$ and the method of this paper provides its norm: $$|t_6(\Po_{-996})| = \frac{2^{16}7^{12}71^{4}83^2}{17^{6}29^{6}41^{6}}.$$
In addition, this method should generalize even further to computing norms of Hauptmoduli on higher genus Shimura curves.



\end{section}

\begin{section}{Shimura Curves} \label{Chapter2}

\begin{subsection}{Quaternion Algebras}

Quaternion algebras have a long history of study so we will only provide a brief summary of the important facts. For a more thourough exploration of quaternion algebras see \cite{CRM}, \cite{Johan}, and \cite{Vig}.

A rational quaternion algebra $B$ is a central simple algebra of dimension $4$ over $\Q$ and is either isomorphic to $\M(\Q)$ or is a skew field. In the latter case, $B$ is called a division algebra. For each prime $p$, $B_p = B \otimes_\Q \Q_p$ is a $\Q_p$-algebra. If $B_p$ is a division algebra, then $B$ is said to be ramified at $p$. If $B_p$ is not a division algebra, then $B_p \simeq \M(\Q_p)$. A quaternion algebra is called definite (indefinite) if it ramifies (is not ramified) at the infinite prime.

The (reduced) discriminant $D=D(B)$ of a quaternion algebra $B$ is given as the product of all finite ramified primes of $B$. Given an even number of finite or infinite primes, there exists a quaternion algebra over $\Q$ ramified exactly at those places. Further, two quaternion algebras are isomorphic if and only if they have the same discriminant.


\begin{prop}[Proposition 3.1 of \cite{Johan}]\label{qDB}
Let $B$ be an indefinite quaternion algebra over $\Q$ with $D = p_1 \cdots p_{2r}$. Choose $q$ to be a prime such that $q \equiv 5 \mod 8$ and $(\frac{q}{p_i})=-1$ for every $p_i>2$. Then $B \simeq \Q(\alpha,\beta)$ where $\alpha\beta=-\beta\alpha$ and $\alpha^2 = q$, $\beta^2=D$. We denote this by $B = \left(\frac{q,D}{\Q}\right)$.
\end{prop}

There are many ways to embed $B=\left(\frac{a,b}{\Q}\right)$ into a matrix algebra over an extension of $\Q$. The one that we use in this paper is $$\phi_b: B \hookrightarrow \M(\Q(\sqrt{b}))$$ given by
\begin{eqnarray*}
\phi_b(\alpha) = \begin{pmatrix}0 & a \\ 1 & 0\end{pmatrix}, & &
\phi_b(\beta)  = \begin{pmatrix}\sqrt{b} & 0 \\ 0 & -\sqrt{b}\end{pmatrix} .
\end{eqnarray*}

There is a natural involution on $x = x_0+x_1\alpha+x_2\beta+x_3\alpha\beta$ given by $$\overline{x} = x_0-x_1\alpha-x_2\beta-x_3\alpha\beta.$$ This involution allows one to define the (reduced) trace and (reduced) norm as
\begin{eqnarray*}
\begin{array}{rcccl}
\tr(x) &=& x + \overline{x} &=& 2x_0, \\
\no(x) &=& x\overline{x} &=& x_0^2-ax_1^2-bx_2^2+abx_3^2.
\end{array}
\end{eqnarray*}
Under the above embedding, these are just the usual matrix trace and determinant. 


\end{subsection} \begin{subsection}{Maximal Orders}

\begin{defn}Let $\K$ be either $\Q$ or $\Q_p$ and $\RK$ its ring of integers. An $\RK$-order $\Oh$ in a quaternion algebra $B$ over $\K$ is an $\RK$-ideal that is a ring. Equivalently, an $\RK$-order $\Oh$ is a ring whose elements have trace and norm in $\RK \subset \Oh$, and $\Oh \otimes_\RK \K = B$. A maximal order is an order that can not be properly contained in another order.
\end{defn}

%
In general, $B$ does not have a unique maximal order. In fact, if $\omega \in B^\times$ and $\Oh$ is a maximal order, then $\omega\Oh\omega^{-1}$ is also a maximal order. However, when $B$ is indefinite, the conjugacy class of maximal orders is unique. 

\begin{prop}[Proposition 3.2 of \cite{Johan}]\label{MaxOrd} For $B$ as in Proposition \ref{qDB} with $\alpha^2 = q$ and $\beta^2=D$, every maximal order is conjugate to
\begin{eqnarray*}
\Oh = \Z + \Z e_1 + \Z e_2 + \Z e_1e_2
\end{eqnarray*}
where
\begin{eqnarray*}
e_1 &=& \frac{1+\alpha}{2}, \\
e_2 &=& \frac{m\alpha+\alpha\beta}{q}, \\
D &\equiv& m^2 \mod q.
\end{eqnarray*}
\end{prop}

When $p$ is a ramified prime, there is a unique maximal order $\Oh_p \subset B_p$ and it is given by
\begin{eqnarray*}
\Oh_p = \{\omega \in B \mid (\ord_p \circ \no)(\omega) \ge 0\}.
\end{eqnarray*}
Hence its group of units is given by
\begin{eqnarray*}
\Oh_p^\times = \{\omega \in B^\times \mid (\ord_p \circ \no)(\omega) = 0\}.
\end{eqnarray*}
Moreover, one can choose a uniformizer $\pi_p \in B_p^\times$ such that $B_p^\times = \Oh_p^\times \rtimes \pi_p^\Z$ with $(\ord_p \circ \no)(\pi_p)=1$ and $\pi_p^2=p$.

Define the normalizer of an order as
\begin{eqnarray*}
\No_{B^\times}(\Oh) = \{\omega \in B^\times \mid \omega\Oh\omega^{-1} \subset \Oh\}.
\end{eqnarray*}
The units of an order $\Oh$ are a subgroup of $\No_{B^\times}(\Oh)$ and are related by the following lemma.

\begin{lemma}[\cite{Vig}]\label{NormDegree} Let $d(B)$ denote the number of ramified primes of $B$. Then
\begin{eqnarray*}
\No_{B^\times}(\Oh)/(\Q^\times\Oh^\times) \simeq (\Z/2\Z)^{d(B)}.
\end{eqnarray*}
\end{lemma}

\end{subsection} \begin{subsection}{Shimura Curves and CM Points}\label{ShCM}

From now on let $B=\left(\frac{q,D}{\Q}\right)$ with $\alpha^2 = q$ and $\beta^2 = D$ as in Proposition \ref{qDB}. Fix the embedding of $B \hookrightarrow \M(\R)$ given by $\phi_{D}$ and the maximal order $\Oh$ as in Proposition \ref{MaxOrd}. Define the following subgroups of $B^\times$,
\begin{eqnarray*}
\Gamma = \Oh^\times, & &
\Gamma^* = \No_{B^\times}(\Oh).
\end{eqnarray*}
Their images under $\phi_{D}$ are discrete subgroups of $B^\times \subset \GL(\R)$, and they act on $\h^{\pm} = \Pro(\C)-\Pro(\R)$ via fractional linear transformations. Define $\X$ and $\Xs$ to be the Shimura curves 
\begin{eqnarray*}
\X = \X_D = \Gamma\bs \h^{\pm}, & & \Xs = \Xs_D = \Gamma^*\bs \h^{\pm}.
\end{eqnarray*}
When $B$ is an indefinite division algebra, $\X$ and $\Xs$ are compact Riemann surfaces with no cusps. Also, Lemma \ref{NormDegree} implies that $\X$ is a covering space of $\Xs$ of degree $2^{d(B)}$.

Fix a quadratic imaginary field $\kay$ such that if $p \mid D$ then $p$ does not split in $\kay$. Then there are many embeddings $\iota: \kay \hookrightarrow B$. However, all of the embeddings are conjugate to each other \cite{Vig}.

\begin{defn}
The image $\iota(\kay^\times) \rightarrow B^\times/\Q^\times \subset \text{PGL}_2(\R)$ has a unique fixed point on $\h^+$. A complex-multiplication (CM) point of $\X$ (resp., $\Xs$) is the $\Gamma$-orbit ($\Gamma^*$-orbit) of such a point. It is said to have discriminant equal to the field discriminant of $\kay$.
\end{defn}

Since all embeddings are conjugate, a CM point is independent of the embedding. In the classical case of $B=\M(\Q)$, the CM points are irrational imaginary solutions to integral quadratic equations with the corresponding discriminant.

\end{subsection} \begin{subsection}{Involutions on $\Xs_D(l)$}

In this section, we summarize the method used in \cite{Elk} to calculate the coordinates of rational CM points on $\Xs$. 
Let $\Po_\Delta$ be the CM point with discriminant $\Delta<0$ and let $R \subset \kay$ be the maximal order in the quadratic imaginary field of discriminant $\Delta$.

\begin{prop}[\cite{Elk}]
$\Po_\Delta$ is a rational point on $\Xs_D$ if and only if the class group of $\kay$ is generated by ideals $I \subset R$ such that $I^2 = (p)$ for some $p \mid D$.
\end{prop}

This implies that for a rational CM point, the class group of $\kay$ is isomorphic to $(\Z/2\Z)^r$ where $r \le d(B)$. In the case of $d(B)=2$, all such fields are known, and thus the rational CM points can be identified. (See Table \ref{Rat6} for $D=6$ and Table \ref{Rat10} for $D=10$.)

Now let $l$ be a prime not dividing $D$, so that $B \otimes_\Q \Q_l \simeq \M(\Q_l)$. Define 
$$\Gamma^*(l) = \{\gamma \in \Gamma^* \mid \gamma \equiv \pm 1 \mod l\}$$
and the congruence subgroup $\Gamma_0^*(l)$ in the same fashion as its classical counterpart. Then the curves
\begin{eqnarray*}
\begin{array}{lcr}
\Xs_D(l) = \Gamma^*(l)\bs\h^\pm, & & \Xs_{D,0}(l) = \Gamma^*_0(l)\bs\h^\pm \end{array}
\end{eqnarray*}
are coverings of $\Xs_D$ whose points are also associated to abelian varieties. From the geometric structure, $\Xs_{D,0}(l)$ inherits an involution $w_l: \Xs_{D,0}(l) \rightarrow \Xs_{D,0}(l)$ which preserves the set of rational CM points.

In the case of $D=6$, the image of $\Gamma^* \hookrightarrow \text{PGL}_2(\R)$ is generated by three elements and is called a triangle group. An area calculation \cite{Elk} shows that $\Xs_6$ has genus $0$. Any coordinate map $t_6:\Xs_6\rightarrow\Pro^1$ is defined up to a $\text{PGL}_2(\R)$ action, so such a map is only well-defined once its values at three points have been given. Since there are three distinguished elements of $\Gamma^*$, the coordinate map is defined to take the values of $0$, $1$, $\infty$ at $\Po_{-4}$, $\Po_{-24}$, $\Po_{-3}$, the CM points associated to the three generators.

The covering curves $\Xs_{6,0}(l)$, for $l=5,7,13$ have genus $0$ and $w_l$ can be expressed explicitly as a rational function. Then by examining the fixed points of $w_l$ and the $w_l$-orbits of $0$, $1$, and $\infty$, Elkies was able to compute the coordinates of 17 of the 27 rational CM points (see Table \ref{Rat6}).

In order to compute the remaining ten rational CM points using this method, involutions on $\Xs_{6,0}(l)$ for higher $l$ are needed. However, these curves have genus greater than $0$ and explicit expressions for $w_l$ are unknown. Instead, Elkies used numerical techniques to calculate the coordinates to an arbitrary precision. He then recognized them as fractional values through continued fractions and their behavior under standard transformations. For example, one expects that the factorizations of both $t_6(\Po_\Delta)$ and $t_6(\Po_\Delta)-1$ should only contain small primes to large powers.

\end{subsection} \end{section} 

\begin{section}{Quadratic Spaces and Lattices} \label{Chapter3}
For a given indefinite quaternion algebra $B$, define the $\Q$-vector space
\begin{eqnarray*}
V = \{ x \in B \mid \tr(x)=0 \}.
\end{eqnarray*}
There is a natural quadratic form on $V$ given by $Q(x)  = \no(x) = -x^2$. Let $(x,y) = \tr(x\overline{y})$ denote the associated inner product 
which has signature $(1,2)$. 

\begin{subsection}{The Lattice $\Oh \cap V$}\label{VLat}
Define the lattice $L = \Oh \cap V$. Let $L^\vee$ be the $\Z$-dual of $L$ 
and consider
$L_p^\vee / L_p$ where $L_p = L \otimes_\Z \Z_p$. 

For $p \nmid D$ and $p$ odd, there is an isomorphism 
$B_p \simeq \M(\Q_p)$ such that $\Oh_p \simeq \M(\Z_p)$. Then  $L_p$ is the
set of trace zero elements of $\M(\Z_p)$ and $L_p^\vee / L_p$ is trivial. Thus
\begin{eqnarray*}
L^\vee/L \simeq \prod\limits_{p \mid 2D} L_p^\vee/L_p.
\end{eqnarray*}

Now consider $p \mid D$ and $p$ odd. Let $\delta \not\in \Z_p^\times$, $\delta^2 \in \Z_p^\times$ and $\Z_{p^2} = \Z_p + \Z_p\delta$ be the ring of integers in the unramified quadratic extension of $\Q_p$ with Galois automorphism $\sigma$. Then 
\begin{eqnarray}
\begin{array}{lcr}
L_p = \Z_p\delta + \Z_p\pi_p + \Z_p\delta\pi_p, & & L_p^\vee = \Z_p \delta + p^{-1}\Z_{p^2}\pi_p. \label{Lp}
\end{array}
\end{eqnarray}
Since $\frac{1}{p}\Z_{p^2}/\Z_{p^2} \simeq \F_{p^2}$, the field of $p^2$ elements, there is an isomorphism
\begin{eqnarray*}
\begin{array}{lcr}
\F_{p^2} \stackrel{\sim}{\rightarrow} L_p^\vee/L_p, & & \tilde{\vpi} \mapsto \vpi\pi_p^{-1} + L_p.
\end{array}
\end{eqnarray*}
Under this isomorphism, the quadratic form $Q$ induces the function
\begin{eqnarray*}
Q(\tilde{\vpi}) = \vpi\vpi^\sigma p^{-1} \mod \Z_p,
\end{eqnarray*}
which is equivalent to the norm map $\no: \F_{p^2} \rightarrow \F_p$ via $\F_p \stackrel{\sim}{\rightarrow}\frac{1}{p}\Z_p/\Z_p$.

The case of $p=2$ has
$L_2^\vee = \frac{1}{2}L_2$. This time the isomorphism is 
\begin{eqnarray}
\begin{array}{lcr}
\F_2 \oplus \F_4 \stackrel{\sim}{\rightarrow} L_2^\vee/L_2, & & (\tilde{w},\tilde{\vpi}) \mapsto w\frac{\sqrt{5}}{2}+\vpi\pi_2^{-1}+L_2, \label{L2}
\end{array}
\end{eqnarray}
and $Q$ induces the function
\begin{eqnarray*}
Q(\tilde{w},\tilde{\vpi}) = -\frac{1}{4}w^2-\frac{1}{2}\no(\vpi) \mod \Z_2.
\end{eqnarray*}
This surjects onto $\frac{1}{4}\Z/\Z$, given by whether or not each of the components is nonzero. 

\begin{prop}\label{DualSize}
Let $D_0$ be the odd part of $D$. Then
\begin{eqnarray*}
|L^\vee/L| = 8D_0^2.
\end{eqnarray*}
\end{prop}

\begin{prop} Let $B_p^\times$ act on $L^\vee_p/L_p$ via conjugation. Then the $B_p^\times$ orbits of $L_p^\vee/L_p$ for odd $p \mid D$ (resp., $p=2$) are indexed by elements of $\F_p$ (\ $\F_4$).\label{QOrb}
\end{prop}
\begin{proof}
For odd $p$, write $B_p^\times$ as
\begin{eqnarray*}
B_p^\times = (\Oh_p^\times \cup \Oh_p^\times\pi_p)p^\Z.
\end{eqnarray*}
First, the powers of $p$ are central and hence act trivially. Then by (\ref{Lp})
\begin{eqnarray*}
L_p^\vee/L_p \stackrel{\sim}{\rightarrow}\Oh_p/\pi_p\Oh_p.
\end{eqnarray*}
Thus the elements of $\Oh_p^\times$ act through their image under the reduction map $\Oh_p \rightarrow \F_{p^2}$. More explicitly, $\tilde{\vpi}\in\F_{p^2}^\times$ acts via left multiplication by $\vpi/\vpi^\sigma$. 
However, this is just the action of $\F_{p^2}^1=\ker(\no: \F_{p^2}^\times \rightarrow \F_p^\times)$. Lastly, $\pi_p$ acts by $\sigma$, and so there is a surjection
\begin{eqnarray*}
B_p^\times \twoheadrightarrow \F_{p^2}^1 \rtimes \langle\sigma\rangle.
\end{eqnarray*}
Hence the orbits of $B_p^\times$ are indexed by the elements of $\F_p$.

For $p=2$, the action of $B_2^\times$ preserves the first component of (\ref{L2}) and acts on the second component the same way it did in the odd $p$ case. So again the orbits are indexed by the four values of $Q$.
\end{proof}

\end{subsection} \begin{subsection}{The Order of the Orbits of $\Gamma^*$}\label{OrderOrbits}
Define the set
\begin{eqnarray*}V(t) = \{x \in V\ |\ Q(x) = t \}\end{eqnarray*}
and $L(t) = L \cap V(t)$.
The discrete groups $\Gamma$ and $\Gamma^*$
both act on $L$ by conjugation, and 
the order of $\Gamma^*$-orbits in $L(t)$ will play an important role in Section \ref{Chapter7}.

Let $0>\Delta \in \Z$ be the field discriminant 
of $\kay = \Q(\sqrt{-t})$, and set $-4t=n^2\Delta$. Then the order $\Z[\sqrt{-t}]$
has discriminant $-4t$. Hence, its conductor is $n$, and if any other order
$R$ in $\kay$ contains $\Z[\sqrt{-t}]$, then the conductor of $R$ divides $n$.

Set
\begin{eqnarray*}\E = \Hom_{\Q\text{-alg}}(\kay, B).\end{eqnarray*}
Assume that for every prime $p \mid D$, $p$ is nonsplit in $\kay$ so that $\E$ is nontrivial.
For every $x \in L(t)$, define
$\iota_x \in \E$ by $\iota_x(\sqrt{-t}) = x$. 
For $\iota \in \E$, $\iota^{-1}(\Oh \cap \iota(\kay))$ is an order
in $\kay$. Let $\cond(\iota)$ denote the conductor of this order and 
define
\begin{eqnarray*}\E(c) = \{ \iota \in \E \mid \cond(\iota)=c\}.\end{eqnarray*}
For $x \in L$, define $\cond(x) = \cond(\iota_x)$ and let
\begin{eqnarray*}L(t,c) = \{x \in L(t) \mid \cond(x) = c \}.\end{eqnarray*}
%
%
Then for a fixed $t$ and $c$, there is a bijection $L(t,c) \stackrel{\sim}{\rightarrow} \E(c)$ given by $x \mapsto \iota_x$ and $\Gamma^*$ acts on $L(t,c)$ via conjugation.
%
This action is compatible with the action on $\E(c)$, therefore
\begin{eqnarray}
\Gamma^* \bs L(t,c) \stackrel{\sim}{\rightarrow} \Gamma^* \bs \E(c).\label{LtoE}\label{bij1}
\end{eqnarray}

To determine the set of $\Gamma^*$-orbits in $L(t,c)$, we examine the right-hand side of (\ref{LtoE}).
Let $R$ be the ring of integers of an imaginary quadratic field $\kay$. Fix an embedding $\iota_0:\kay \hookrightarrow B$ with $\cond(\iota_0)=1$,
i.e. $\iota_0(R) \subset \Oh$.
Since all embeddings of $\kay$ into $B$ are conjugate, there is a bijection
\begin{eqnarray*}
B^\times /\kay^\times  \stackrel{\sim}{\rightarrow} \E, \\
\omega \mapsto Ad(\omega) \circ \iota_0.
\end{eqnarray*}
Then
\begin{eqnarray*}\Gamma^* \bs B^\times /\kay^\times  \stackrel{\sim}{\rightarrow} \Gamma^* \bs \E,\end{eqnarray*}
where the action of $\Gamma^*$ on $B^\times /\kay^\times $ is left multiplication. Define
\begin{eqnarray*}B^\times (c) = \{\omega \in B^\times  \mid \cond(Ad(\omega) \circ \iota_0)=c \}\end{eqnarray*}
so that
\begin{eqnarray}\Gamma^* \bs B^\times (c)/\kay^\times  \stackrel{\sim}{\rightarrow} \Gamma^* \bs \E(c).\label{bij2}\end{eqnarray}

Let $\Ord = \Ord(B)$ be the set of all maximal orders of $B$.
For any $\Oh \in \Ord$, define the conductor of $\Oh$ to be the conductor of
$\iota_0^{-1}(\Oh \cap \iota_0(\kay))$. Define for $\omega \in B^\times $, $\Oh_\omega = \omega^{-1}\Oh \omega \in \Ord$. Then
the conductor of $\Oh_\omega$ is $\cond(\omega)$.

The action of $B^\times_{\Af} = (B \otimes_\Q \Af)^\times$ on $\Ord$ via
\begin{eqnarray*}\xi \cdot \Oh_\omega = \xi^{-1}\widehat{\Oh_\omega}\xi \cap B\end{eqnarray*}
where $\widehat{\Oh} = \Oh \otimes_\Z \widehat{\Z}$.
is transitive, thus
\begin{eqnarray*}
\No_{B^\times_{\Af} }(\widehat{\Oh}) \bs B^\times_{\Af}  &\stackrel{\sim}{\rightarrow}& \Ord, \\
\xi &\mapsto& \xi^{-1}\widehat{\Oh} \xi \cap B. \end{eqnarray*}
Furthermore, the double cosets
\begin{eqnarray*}\No_{B^\times_{\Af} }(\widehat{\Oh}) \bs B^\times_{\Af}  / B^\times \end{eqnarray*}
correspond to the $B^\times $-conjugacy classes of the maximal orders in $B$.
Since $B$ is an indefinite quaternion algebra, all maximal orders of $B$ are conjugate. Thus
\begin{eqnarray*}\No_{B^\times_{\Af} }(\widehat{\Oh}) \bs B^\times_{\Af} \simeq \No_{B^\times }(\Oh) \bs B^\times.  
\end{eqnarray*}
Let $\Ord(c) \subset \Ord$ be the subset of orders with conductor $c$. Then, with notations as before, 
\begin{eqnarray}
\No_{B^\times }(\Oh) \bs B^\times (c) \stackrel{\sim}{\rightarrow} \Ord(c), \label{bij3}
\end{eqnarray}
and the $\kay^\times_{\Af}= (\kay \otimes_\Q \Af)^\times$ action on $\Ord$ given by $\xi \cdot \Oh = \xi^{-1}\widehat{\Oh}\xi \cap B$ preserves $\Ord(c)$.

From the Chevalley-Hasse-Noether theorem, for a given $\Oh_c \in \Ord(c)$ there is a bijection
\begin{eqnarray}
\No_{B^\times_{\Af} }(\Oh_c) \cap \kay^\times_{\Af}  \bs \kay^\times_{\Af}  \stackrel{\sim}{\rightarrow} \Ord(c) \label{bij4}
\end{eqnarray}
given by the orbit of $\Oh_c$ under the transitive action of $\kay^\times_{\Af} $. Then the composition of the bijections in (\ref{bij1}), (\ref{bij2}), (\ref{bij3}), and (\ref{bij4}) yield 
\begin{eqnarray*}
\Gamma^* \bs L(t,c) \stackrel{\sim}{\leftrightarrow} \No_{B^\times_{\Af} }(\Oh_c) \cap \kay^\times_{\Af}  \bs \kay^\times_{\Af} / \kay^\times.
\end{eqnarray*}

Let $\Delta_0$ be the product of all the primes that ramify in $\kay$ and define
\begin{eqnarray*}
\delta(\Delta_0,D)=\#\{p \text{ prime}\ |\ p\mid\gcd(\Delta_0,D)\} -\left\{\begin{array}{ll} 1 & \text{if } \Delta_0 \mid D \\
0 & \text{otherwise}\end{array}\right.
\end{eqnarray*}
\begin{theorem} Let $R_c \in \kay$ be the order of conductor $c$, then
$$[\widehat{R_c}^\times \bs \kay^\times_{\Af}  / \kay^\times : \No_{B^\times_{\Af} }(\Oh_c) \cap \kay^\times_{\Af}  \bs \kay^\times_{\Af} / \kay^\times] = 2^{\delta(\Delta_0,D)}.$$
\end{theorem}
\begin{proof}
For a prime $p \nmid D$, $$\No_{B^\times_{p} }(\Oh_c) = \Oh_{c,p}^\times\Q_p^\times$$ thus 
\begin{eqnarray}
\No_{B^\times_{p} }(\Oh_c) \cap \kay_p^\times = R_{c,p}^\times\Q_p^\times. \label{NkB}
\end{eqnarray} 

For primes $p \mid D$, $\No_{B^\times_{p} }(\Oh_c) = B^\times_{p}$. When $p$ is inert in $\kay$, (\ref{NkB}) still holds. However, when $p$ is ramified in $\kay$, $$\No_{B^\times_{p} }(\Oh_c) \cap \kay_p^\times = R_{1,p}^\times\Q_p^\times \cup R_{1,p}^\times\Q_p^\times \pi_p$$ where $\pi_p^2=p$. 

Altogether, then, there is a surjection
\begin{eqnarray*}
\widehat{R_c}^\times \bs \kay^\times_{\Af}  / \kay^\times \twoheadrightarrow \No_{B^\times_{\Af} }(\Oh_c) \cap \kay^\times_{\Af}  \bs \kay^\times_{\Af} / \kay^\times
\end{eqnarray*}
given by modding out by the subgroup generated by the elements $(1,...,1,\pi_p,1,...)$ for $p$ ramified in both $B$ and $\kay$. The size of this subgroup is $2^{\delta(\Delta_0,D)}$.
\end{proof}

\begin{cor}\label{A2}Let $h(c^2\Delta)$ be the ideal class number of the order of conductor $c$ in the quadratic field of discriminant $\Delta$ and $-4t=n^2\Delta$ as before. Then
\begin{eqnarray*}
|\Gamma^* \bs L(t)| = 2^{-\delta(\Delta_0,D)}\sum\limits_{c \mid n} h(c^2\Delta).
\end{eqnarray*}
where $h(c^2\Delta)$ is the class number of $R_c$, the order of conductor $c$ in $\kay$.\end{cor}
\begin{proof}
This follows from recognizing $\widehat{R_c}^\times \bs \kay^\times_{\Af}  / \kay^\times$ as the desired ideal class group and noting that $L(t) = \coprod\limits_{c \mid n} L(t,c)$.
\end{proof}

\end{subsection} \end{section} \begin{section}{Borcherds Forms} \label{Chapter4}

\begin{subsection}{Rational Quadratic Divisors}
Let $\D$ be the space of oriented negative $2$-planes in $V$. Call $[z_1,z_2] \in \D$ a proper basis if $(z_1,z_1)=(z_2,z_2)=-1$ and $(z_1,z_2)=0$. 
In addition, define 
\begin{eqnarray*}
{\cal Q}=\{v \in V(\C) \mid (v,v)=0,(v,\overline{v})<0\}/\C^\times.
\end{eqnarray*}
This is an open subset of a quadric in $\Pro(V(\C))$.
Recall that $B=\left(\frac{q,D}{\Q}\right)$ with $\alpha^2=q$ and $\beta^2=D$ and 
let $V$ have the canonical basis $\{\alpha, \beta, \alpha\beta\}$. Then there is a pair of bijections
\begin{eqnarray*}
\h^{\pm} \stackrel{w}{\longrightarrow} {\cal Q} \stackrel{\sigma}{\longleftarrow} \D(\R)
\end{eqnarray*}
%
where the maps are given by
\begin{eqnarray}
w(z) &=& \left(\frac{q-z^2}{2q}\right)\alpha + \left(\frac{z}{\sqrt{D}}\right)\beta + \left(\frac{q+z^2}{2q\sqrt{D}}\right)\alpha\beta, \label{HalfPlaneMap}\\
\sigma([z_1,z_2]) &=& z_1 - i z_2.
\end{eqnarray}
%


Write $\D = \D^+ \cup \D^- $ where $\D^+$ (resp., $\D^-$) are the planes with positive (negative) orientation. For $x \in V(\Q)$ define
\begin{eqnarray*}
\D_x = \{z \in \h^{\pm} \mid (x,w(z))=0\}.
\end{eqnarray*}
By (\ref{HalfPlaneMap}), for $x = x_1 \alpha + x_2 \beta + x_3 \alpha\beta$,
\begin{eqnarray}
(x,w(z)) = \left(\frac{x_1+x_3\sqrt{D}}{2}\right)z^2 - (x_2 \sqrt{D})z - \frac{q(x_1-x_3\sqrt{D})}{2}. \label{DxEq}
\end{eqnarray}
Hence
\begin{eqnarray*}
\D_x = \left\{\frac{x_2 \sqrt{D} \pm \sqrt{-Q(x)}}{x_1+x_3\sqrt{D}} \right\}.
\end{eqnarray*}
Let $\D_x^\pm = \D_x \cap \D^\pm$.

\begin{prop}\label{DxFixed}
For $x \in V$ with $Q(x)>0$, $\D_x$ is the set of fixed points of the image of $x$ in $\text{PGL}_2(\R)$ under the embedding $\phi_{D}$.
\end{prop}
\begin{proof}
Let $x=x_1 \alpha + x_2 \beta + x_3 \alpha\beta$. Then
\begin{eqnarray*}
\phi_{D}(x) = \begin{pmatrix}x_2 \sqrt{D} & q(x_1-x_3\sqrt{D}) \\ x_1+x_3\sqrt{D} & -x_2 \sqrt{D}\end{pmatrix}.
\end{eqnarray*}
A fixed point, $z$, of this matrix satisfies
\begin{eqnarray*}
zx_2 \sqrt{D} + q(x_1-x_3\sqrt{D}) = z^2 (x_1+x_3\sqrt{D})  - zx_2 \sqrt{D}.
\end{eqnarray*}
This is equivalent to (\ref{DxEq}).
\end{proof}

\begin{defn}\label{ZDef}
Let $G=\Gamma$ or $\Gamma^*$ and let $G\eta$ denote the $G$-orbit of $\eta \in L^\vee/L$.
The rational quadratic divisor $Z(d,\eta;G)$ is given by
\begin{eqnarray*}
Z(d,\eta;G) = \mathop{\mathop{\sum_{x \in L^\vee \cap V(d)}}_{x+L \in G\eta}}_{\text{{\rm mod }} G}\text{{\rm pr}}_G(\D_x^+),
\end{eqnarray*}
where $\text{{\rm pr}}_G:\D^+ \rightarrow G \bs \D^+$ and each point is counted with weight $|\text{Stab}(x)|^{-1}$.\end{defn}
\noindent For more details on this definition in the case of $G=\Gamma$, see the Appendix of \cite{SSKMSRI}. 

\end{subsection} \begin{subsection}{Borcherds Forms}

Let $H=\text{GSpin}(V)$. Viewed as an algebraic group, $H(\mathcal{A}) \simeq (B \otimes_\Q \mathcal{A})^\times$ for any $\Q$-algebra $\mathcal{A}$. Let $K \subset H(\Af)$ be a compact open set such that $H(\A) = H(\Q)H(\R)^+K$ where $H(\R)^+$ is the component of $H(\R)$ that contains the identity. 
\begin{defn}\label{BorForm}
A modular form of weight $k \in \Z$ on $\D \times H(\Af)/K$ is a function $\Psi:\D\times H(\Af)/K \rightarrow \C$ such that
$$\Psi(\gamma z, \gamma h) = j(\gamma,z)^k \Psi(z,h)$$ for all $\gamma \in H(\Q)$, where $j(\gamma,z)$ is the automorphy factor given 
in \cite{SSKBor}.
\end{defn}
\noindent The cases we will focus on have $k=0$ and thus the automorphy factor will be inconsequential.

Let $L$ be a lattice and $F$ be a modular form valued in $\C[L^\vee/L]$ with Fourier expansion given by
\begin{eqnarray}
F(\tau) = \sum\limits_{\eta \in L^\vee/L} \sum\limits_{m \in \Q} c_\eta(m)\q^m e_\eta \label{Fmod}
\end{eqnarray}
where $\{e_\eta\}_{\eta\in L^\vee/L}$ form the basis of $\C[L^\vee/L]$. Since $\Gamma$ and $\Gamma^*$ act on $L^\vee/L$, they also act via linearity on the algebra $\C[L^\vee/L]$ and the function $F$.

\begin{defn}
For a lattice $L$ with signature $(n,2)$, a Borcherds form $\Psi(F)$ is a meromorphic modular form on $\D \times H(\Af)/K$ arising from the regularized theta lift of a weight $1-\frac{n}{2}$ meromorphic modular form $F$ as in (\ref{Fmod}) with $c_\eta(m) \in \Z$ for $m \le 0$. See \cite{JS}, \cite{SSKBor}, \cite{BorGrass}.
\end{defn}

Borcherds forms have the following key properties. 

\begin{theorem}[Theorem 1.3 of \cite{SSKBor}]\label{DivForm} Assume $F$ is given as in (\ref{Fmod}) and is $\Gamma^*$ invariant.
\begin{enumerate}
\item[1)] The weight of $\Psi(F)$ is $c_0(0)$.
\item[2)] $\text{{\rm div}}(\Psi(F)^2)=\sum\limits_{\eta \in L^\vee/L} \sum\limits_{m>0} c_\eta(-m)Z(m,\eta;\Gamma^*).$
\end{enumerate}
\end{theorem}

\end{subsection} \begin{subsection}{Adelic View}

We can rephrase some of the definitions from Section \ref{ShCM} from an adelic point of view. This will allow the machinary of Borcherds forms to apply to the computation of singular moduli on $\X_D$ and $\Xs_D$.

Let $K_\Gamma$ be the compact open set $\widehat{\Oh}^\times \subset H(\Af)$. Then $\Gamma = H(\Q)\cap H(\R)^+K_\Gamma$.
Let $K_{\Gamma^*}$ be defined analagously. Then $\X_D$ and $\Xs_D$ are given by
\begin{eqnarray*}
\X_D \simeq \Gamma\bs \D \simeq H(\Q)\bs(\D \times H(\Af)/K_\Gamma), \\
\Xs_D \simeq \Gamma^*\bs \D \simeq H(\Q)\bs(\D \times H(\Af)/K_{\Gamma^*}).
\end{eqnarray*}
\noindent Notice that $\X_D$ and $\Xs_D$ are natural domains for weight-$0$ Borcherds forms.

The CM points can be viewed adelically as well. An element $x\in V(\Q)$ with positive norm gives rise to the decomposition of $V$ as
$V = \Q x \oplus U$
where $U = x^\perp$ is a negative plane. This splitting corresponds to a two-point set $\D_x$. As a rational inner product space $U \simeq \kay$ for some quadratic imaginary field $\kay$ with quadratic form given by a constant times the norm on $\kay$. Set $T\simeq\text{GSpin}(U)$. Then, with $\iota_x$ as in  Section \ref{OrderOrbits}, $T(\Q) \simeq \iota_x(\kay^\times) \subset H(\Q)$ and the CM points are the image of 
\begin{eqnarray}
Z_{\Gamma^*}(U) = T(\Q)\bs(\D_x \times T(\Af)/K_{\Gamma^*}) \hookrightarrow \Xs_D. \label{CMZ}
\end{eqnarray}
The degree of this $0$-cycle is given in Chapter 3 of \cite{KRYBook} as
\begin{eqnarray*}
|Z_{\Gamma^*}(U)| = 2 \cdot \sum_{c \mid n} \frac{h(c^2\Delta)}{w(c^2\Delta)} \cdot \prod_{p|D}(1-\chi_\Delta(p)) 
\end{eqnarray*}
where $w(c^2\Delta)$ is the number of units in $R_c$ and $\chi_\Delta$ is the associated Dirichlet character for $\kay$ given by the Kronecker symbol, $\chi_\Delta(n) = \left(\frac{\Delta}{n}\right)$. 
\end{subsection} \begin{subsection}{Borcherds Forms at CM Points}

Recall that $L=\Oh \cap V$ is a lattice in $V$ corresponding to a fixed maximal order $\Oh$. Then there are sublattices
\begin{eqnarray*}
\begin{array}{ccc}
L_+ = \Q x \cap L, & & L_- = U \cap L.
\end{array}
\end{eqnarray*}
In general, $L \ne L_-+L_+$, and 
$$L_-+L_+ \subseteq L \subseteq L^\vee \subseteq L_-^\vee+L_+^\vee.$$
Hence an element $\eta \in L^\vee$ decomposes as $\eta= \eta_- + \eta_+$ for $\eta_\pm \in L_\pm^\vee$.


\begin{defn}[\cite{JS}]\label{JSdef} For $\mu \in L_-^\vee/L_-$ and $\psi_\mu = \text{char}(\mu+L_-)$, let $E(\tau,s;\psi_\mu,+1)$ be the incoherent Eisenstein series of weight $1$ with Fourier expansion
\begin{eqnarray*}
E(\tau,s;\psi_\mu,+1) = \sum\limits_{m}A_\mu(s,m,v)\q^m
\end{eqnarray*}
where the Fourier coefficients have Laurent expansions
\begin{eqnarray*}
A_\mu(s,m,v) = b_\mu(m,v)s + O(s^2)
\end{eqnarray*}
at $s=0$. Then for $\eta \in L^\vee/L$ and $m \in \Q$ define
\begin{eqnarray}\label{kappasum}
\kappa_\eta(m) = \sum\limits_{\lambda \in L/(L_+ + L_-)} \sum\limits_{x \in \eta_+ + \lambda_+ + L_+} \kappa^-_{\eta_-+\lambda_-}(m-Q(x)) \label{kappaeq}
\end{eqnarray}
where
\begin{eqnarray}
\kappa^-_\mu(m') &=& \left\{\begin{array}{ll}\lim_{v\rightarrow \infty}b_\mu(m',v) & \text{if }m'>0 \\
k_0(0)\psi_\mu(0) & \text{if }m'=0 \\
0 & \text{if }m'<0\end{array} \right. , \label{kapmindef}\\
k_0(0) &=& \log(|\Delta|)+2\frac{\Lambda'(1+\chi_\Delta)}{\Lambda(1,\chi_\Delta)},
\end{eqnarray}
and $\Lambda(s,\chi_\Delta)$ is the normalized $L$-series $\pi^{-\frac{s+1}{2}}\Gamma\left(\frac{s+1}{2}\right)L(s,\chi_\Delta)$. 
\end{defn}

\begin{theorem}[Corollary 3.4 of \cite{JS}]\label{(1,2)} Assume $c_\eta(m) \in \Z$ for $m \le 0$, $c_0(0)=0$, and that the $0$-cycle $Z_{\Gamma^*}(U)$ defined in (\ref{CMZ}) does not meet the divisor of $\Psi(F)$. Then
\begin{eqnarray}
\frac{1}{|Z_{\Gamma^*}(U)|}\sum\limits_{z\in Z_{\Gamma^*}(U)}\log||\Psi(z,f)||^2 = \frac{-1}{2^{d(B)}} \sum\limits_{\eta}\sum\limits_{m \ge 0}c_\eta(m)\kappa_\eta(m) \label{(1,2)eq}
\end{eqnarray}
where $h(\kay)$ is the ideal class number of the quadratic field $\kay \simeq U$.
\end{theorem}

The power of this theorem lies in the explicit formulas for the right-hand side of (\ref{(1,2)eq}). In Section \ref{Chapter7} we will use this theorem to compute the norms of singular moduli. However, first a supply of appropriate vector-valued modular forms $F$ is needed to serve as the input to the Borcherds construction of $\Psi(F)$.

\end{subsection} \end{section} \begin{section}{Input Forms} \label{Chapter5}

This section is presented in general terms and follows \cite{BorGrass} and \cite{BorRef}. However, rather than appearing redundant, the notation implies how the general theory applies to the  set-up in Sections \ref{Chapter2} through \ref{Chapter4}.

\begin{subsection}{$\MSL(\Z)$ and the Weil Representation}

The Lie group $\SL(\R)$ has a double cover $\MSL(\R)$ with elements of the form
$$\left(\begin{pmatrix} a & b \\ c & d \end{pmatrix}, \pm \sqrt{c \tau + d}  \right).$$
The group structure is given by
$$(G_1, j_1(\cdot))(G_2, j_2(\cdot)) = (G_1G_2, j_1(G_2(\cdot))j_2(\cdot)).$$
The group $\MSL(\Z)$ is defined as the inverse image in $\MSL(\R)$ of $\SL(\Z)$ and
is generated by the two elements
\begin{eqnarray*}
S = \left( \begin{pmatrix} 0 & -1 \\ 1 & 0 \end{pmatrix}, \sqrt{\tau} \right), & & T = \left( \begin{pmatrix} 1 & 1 \\ 0 & 1 \end{pmatrix}, 1 \right),
\end{eqnarray*}
which satisfy
\begin{eqnarray*}
Z = S^2 = (ST)^3 = \left( \begin{pmatrix} -1 & 0 \\ 0 & -1 \end{pmatrix}, i \right).
\end{eqnarray*}
The element $Z$ generates the center of $\MSL(\Z)$ and the quotient by $Z^2$ is $\SL(\Z)$. Also, $\MSL(\Z)$ acts on $\h^\pm$ via its image in $\SL(\Z)$.
Throughout the following, let
\begin{eqnarray}
\gamma = \gamma^\pm = \left(\begin{pmatrix}a & b \\ c & d\end{pmatrix}, \pm \sqrt{c \tau + d} \right) \in \MSL(\Z). \label{gamdef}
\end{eqnarray}

Let $L$ be a lattice with quadratic form $Q'$ and let $L^\vee$ be the dual lattice under the associated inner product. To ease notation, let $\Lambda_L = L^\vee/L$.
Then Milgram's formula gives $\sign(L)$, the signature mod 8 of $L$, via
\begin{eqnarray*}
\sum\limits_{\eta \in \Lambda_L}\expo(Q'(\eta)) = \sqrt{|\Lambda_L|}\expo(\sign(L)/8)
\end{eqnarray*}
where $\expo(x) = e^{2\pi i x}$. 
For $\eta \in \Lambda_L$, let $e_\eta$ denote the corresponding basis element in the group ring $\C [\Lambda_L]$. In \cite{BorGrass},
Borcherds defines the Weil representation $\Weilb$ on the generators of $\MSL(\Z)$ in terms of $Q'$.
However, we will use the dual representation $\Weil = \Weilb^\vee$ since the quadratic form in Sections \ref{Chapter2} through \ref{Chapter4} is actually given by $Q(x) = -Q'(x)$. On the generators $\Weil$ is given by
\begin{eqnarray*}
\Weil(T)e_\eta &=& \expo(-Q(\eta))e_\eta, \\
\Weil(S)e_\eta &=& C_L \sum\limits_{\delta \in \Lambda_L} \expo(-(\eta, \delta))e_\delta
\end{eqnarray*}
where 
\begin{eqnarray*}
C_L = \frac{\expo(\sign(L)/8)}{\sqrt{|\Lambda_L|}} = \frac{1}{|\Lambda_L|}\sum\limits_{\eta \in \Lambda_L}\expo(Q(\eta)).
\end{eqnarray*}
(This approach follows \cite{JS} and \cite{SSKBor}. However most of the results in this section are the dualized versions of those found in \cite{BorRef}.)
Define the level of $\Lambda_L$ to be the smallest integer $N$ such that $NQ(\eta) \in \Z$ for all $\eta \in L^\vee$.
Then the representation $\Weil$ factors through $\MSL(\Z/N\Z)$, the double cover of $\SL(\Z/N\Z)$.
Define the congruence subgroup $\Go(N) \subset \SL(\Z)$ as the preimage of the upper triangular matrices in $\SL(\Z/N\Z)$ and $\MGo(N)$ as its inverse image in $\MSL(\Z)$.
\begin{defn}[\cite{BorRef}]For $\gamma \in \MGo(N)$ define
\begin{eqnarray}
\chi_n(\gamma) &=& \left(\frac{d}{n}\right), \\
\chi_\theta(\gamma^\pm) &=& \left\{\begin{array}{ll}
\pm \left(\frac{c}{d}\right) & d \equiv 1 \mod 4 \\
\mp i\left(\frac{c}{d}\right) & d \equiv 3 \mod 4 \\
\end{array} \right., \label{chitheta} \\
\chi_L(\gamma) &=& \left\{ \begin{array}{ll}
\left(\chi_\theta^{-\rm{sign}(L) + \left(\frac{-1}{|\Lambda_L|}\right)-1}\chi_{|\Lambda_L|2^{\rm{sign}(L)}}\right)(\gamma) & 4 \mid N \\
\chi_{|\Lambda_L|}(\gamma) & 4 \nmid N
\end{array} \right. .\label{Kyle}
\end{eqnarray}
\end{defn}

\begin{theorem}[Theorem 5.4 of \cite{BorRef}]Suppose $\Lambda_L$ has level $N$. If $b$ and $c$ are divisible by $N$ then 
$\gamma \in \MSL(\Z)$ acts on $\C[\Lambda_L]$ by
$$\Weil(\gamma) e_\eta = \chi_L(\gamma)e_{a\eta}.$$
\end{theorem}

\begin{cor}Suppose $\Lambda_L$ has level $N$ and that $\eta \in \Lambda_L$ has norm $0$. Then $\gamma \in \MGo(N)$ acts on the element $e_\eta$ by
$$\Weil(\gamma) e_\eta = \chi_L(\gamma)e_{a\eta}.$$
\end{cor}
\begin{proof}
Any element $\gamma \in \MGo(N)$ can be written as
$$\gamma=T^n \left(\begin{pmatrix} a' & b' \\ c & d \end{pmatrix},\pm \sqrt{c \tau + d}  \right)$$
where $N$ divides $c$ and $b'$. Then $\chi_L$ is trivial on $T$. Since $a' \equiv a \mod N$ and the order of $\eta$
divides $N$, $a'\eta = a \eta$.
\end{proof}

\end{subsection} \begin{subsection}{Vector-Valued Modular Forms}
Define the slash operator of weight $k$ for an element $\gamma \in \MSL(\Z)$ by
$$f|^k_{\gamma^\pm}(\tau) = (\pm \sqrt{c \tau + d})^{2k}f(\gamma \tau).$$

\begin{defn}
Suppose $\rho$ is a representation of $\Gamma \subset \MSL(\Z)$ on a finite dimensional complex vector space $\mathcal{V}$. Then  
$F: \h^\pm \rightarrow \mathcal{V}$ is a vector-valued modular form on $\Gamma$ of weight $k \in \frac{1}{2}\Z$ and type $\rho$ if it is meromorphic
and satsifies
\begin{eqnarray*}
F(\gamma^\pm\tau)=(\pm \sqrt{c \tau + d})^{2k} \rho(\gamma^\pm)F(\tau) 
\end{eqnarray*}
for all $\gamma \in \Gamma$.
\end{defn}

\begin{defn}Suppose $f$ is a scalar-valued weight $k$ modular form on $\MGo(N)$ with character $\chi_L$. Then define a weight $k$ modular form $F_f(\tau)$ valued in $\C[L^\vee/L]$ via
\begin{eqnarray}
F_f(\tau) = \sum\limits_{\gamma \in \MGo(N) \bs \MSL(\Z)}f|^k_\gamma(\tau)\Weil(\gamma^{-1})e_0. \label{Fdef}
\end{eqnarray}
\end{defn}
It can be shown \cite{Bar} that $F_f(\tau)$ is well-defined and is a modular form of type $\Weil$ and weight $k$ on $\MSL(\Z)$. 

\begin{prop}
Let $F_f$ have Fourier expansion as in (\ref{Fmod}). If $m+Q(\eta) \not\in \Z$, then $c_\eta(m) = 0$.
\end{prop}
\begin{proof}
Since $F_f$ is a modular form, 
\begin{eqnarray*}
F_f(\tau+1) &=& \Weil(T)F_f(\tau) \\
\sum\limits_{\eta \in \Lambda_L}\sum\limits_{m\in\Q}c_\eta(m)\q^m\expo(m)e_\eta &=& \sum\limits_{\eta \in \Lambda_L}\sum\limits_{m\in\Q}c_\eta(m)\q^m\Weil(T)e_\eta \\
&=& \sum\limits_{\eta \in \Lambda_L}\sum\limits_{m\in\Q}c_\eta(m)\q^m\expo(-Q(\eta))e_\eta.
\end{eqnarray*}
Thus $m+Q(\eta) \not \in \Z$ implies $c_\eta(m) = 0$.
\end{proof}

\begin{prop}
If $f$ has no poles at finite cusps, then, for $F_f$ as in (\ref{Fmod}), $c_\eta(m) = 0$ for $m<0$ and $\eta \ne 0$.
\end{prop}
\begin{proof}
If $f$ does not have a pole at a finite cusp, then the coordinate function $f|_\gamma^k$ in (\ref{Fdef}) can have a pole only when $\gamma(\infty)=\infty$. However, this is satisfied only by the trivial coset representative
which has $\Weil(\gamma^{-1})e_0=e_0$.
\end{proof}

Now define, as in \cite{BorRef}, $\Lambda_{L,n}$ to be the set of $n$-torsion points and define $\Lambda_L^n$ via the exact sequence
$$0 \rightarrow \Lambda_{L,n} \rightarrow \Lambda_n \rightarrow \Lambda_L^n \rightarrow 0,$$
and
$$\Lambda_L^{n*} = \{\delta \in \Lambda_L^n \mid \begin{array}{cc}(\delta, \eta) = -nQ(\eta) & \forall \eta \in \Lambda_{L,n}\end{array} \}.$$
\begin{lemma}\label{normcond}
For a fixed $n$, either $\Lambda_L^{n*} = \varnothing $ or the membership of $\delta$ into $\Lambda_L^{n*}$ is completely determined by $Q(\delta)$. 
\end{lemma}
\begin{proof}
It suffices to examine the criteria locally at the primes that divide the level $N$. Recall from Section \ref{VLat} that for an odd prime $p$, $\Lambda_{L,p} \simeq \F_{p^2}$ and $Q:\Lambda_{L,p} \rightarrow (1/p)\Z/\Z$. If $p \mid n$, then $(\Lambda_{L,p})_n = \Lambda_{L,p}$ and $(\Lambda_{L,p})^n = \{0\}$. Since $nQ(\delta)=0=(0,\delta)$ for all $\delta \in (\Lambda_{L,p})_n$, then $(\Lambda_{L,p})^{n*} = \{0\}$. If $p \nmid n$, then $(\Lambda_{L,p})_n = \{0\}$ and $(\Lambda_{L,p})^n = \Lambda_{L,p}$. Since $nQ(0)=0=(\delta,0)$ for all $\delta \in (\Lambda_{L,p})^n = \Lambda_{L,p}$, then $(\Lambda_{L,p})^{n*} = \Lambda_{L,p}$. So for odd $p \mid N$,
\begin{eqnarray}\label{oddp}
(\Lambda_{L,p})^{n*} = \left\{\begin{array}{ll} \{\delta \mid Q(\delta)\in(1/p)\Z_p/\Z_p\} & p\nmid n \\ \{0\} & p \mid n \end{array}\right. .
\end{eqnarray}

Now consider $p=2$ where $\Lambda_{L,2} \simeq \F_4 \oplus \F_2$ and $Q:\Lambda_{L,2} \rightarrow (1/4)\Z/\Z$. Suppose $2 \nmid n$. Then $(\Lambda_{L,2})_n = \{0\}$, and $(\Lambda_{L,2})^n = \Lambda_{L,2}$. Since $nQ(0)=0=(\delta,0)$ for all $\delta \in (\Lambda_{L,2})^n = \Lambda_{L,2}$, then $(\Lambda_{L,2})^{n*} = \Lambda_{L,2}$. Now suppose $2 \mid n$. Then $(\Lambda_{L,2})_n = \Lambda_{L,2}$, and $(\Lambda_{L,2})^n = \{0\}$. However, $nQ(\delta)=0=(0,\delta)$ for all $\delta \in (\Lambda_{L,p})$ only when $4 \mid n$. Thus 
\begin{eqnarray}\label{evenp}
(\Lambda_{L,2})^{n*} = \left\{\begin{array}{ll} \{\delta \mid Q(\delta)\in(1/4)\Z_2/\Z_2\} & 2\nmid n \\ \varnothing & 2\ \|\ n \\ \{0\} & 4 \mid n \end{array}\right. .
\end{eqnarray}
Combining (\ref{oddp}) and (\ref{evenp}) into one global statement yields
\begin{eqnarray*}
\Lambda_L^{n*} = \left\{\begin{array}{ll} \{\delta \mid Q(\delta)\in\left(\frac{\gcd(n,N)}{N}\right)\Z/\Z\} & 2\not\!\|\ n \\ \varnothing & 2\ \|\ n \end{array}\right. .
\end{eqnarray*}
Thus the membership of an element is determined by its image under $Q$.
\end{proof}

\begin{lemma}[Lemma 3.1 of \cite{BorRef}] The sum
$$\mathcal{S}_n(\delta) = \sum\limits_{\eta \in \Lambda_L}\expo(-(\eta,\delta)-nQ(\eta))$$
is equal to $0$ when $\delta \not\in \Lambda_L^{n*}$ and has magnitude $\sqrt{|\Lambda_L||\Lambda_{L,n}|}$ otherwise.
\end{lemma}
\begin{lemma}[Lemma 3.2 of \cite{BorRef}] For $\gamma \in \MSL(\Z)$ as in (\ref{gamdef}), $\Weil(\gamma)e_0$ is a linear combination of the elements $e_\delta$ for $\delta \in \Lambda_L^{c*}$.
\end{lemma}
\begin{proof}
Since the coset representatives of $\Gamma_0(N) \bs \SL(\Z)$ can all be chosen of the form $S^{-1}T^{-n}S^{-1}T^{-m}$, it is sufficient to prove this for $\gamma$ of the form $T^mST^nS$ for some $m,n \in \Z$ with $(N,n) = (N,c)$ since any $\gamma$ is a product of an element of this form with an element of $\MGo(N)$ on the right, but $e_0$ is an eigenvector for $\MGo(N)$. Then
\begin{eqnarray}
\Weil(S)e_0 &=& C_L \sum\limits_{\delta \in \Lambda_L} e_\delta, \nonumber \\
\Weil(T^nS)e_0 &=& C_L \sum\limits_{\delta \in \Lambda_L} \expo(-nQ(\delta))e_\delta, \nonumber \\
\Weil(ST^nS)e_0 &=& C_L^2 \sum\limits_{\delta \in \Lambda_L} \sum\limits_{\delta' \in \Lambda_L}\expo(-nQ(\delta)-(\delta', \delta))e_{\delta'}\nonumber \\
&=& C_L^2 \sum\limits_{\delta \in \Lambda_L^{n*}}\mathcal{S}_n(\delta)e_{\delta},\nonumber \\
\Weil(T^mST^nS)e_0 &=& C_L^2\sum\limits_{\delta \in \Lambda_L^{n*}}\mathcal{S}_n(\delta)\expo(-mQ(\delta))e_\delta. \label{pTSTS}
\end{eqnarray}
\end{proof}

\begin{theorem}
If $Q(\delta)=Q(\delta')$, then the $e_\delta$ and $e_{\delta'}$ components of $F_f$ are equal.
\end{theorem}
\begin{proof} This follows from the fact that the coefficient $\mathcal{S}_n(\delta)\expo(-mQ(\delta))$ in (\ref{pTSTS}) depends only on $Q(\delta)$ which, by Proposition \ref{normcond}, is the same for all $\delta \in \Lambda_L^{n*}$.
\end{proof}

\begin{cor}
The modular form $F_f$ is $\Gamma^*$ invariant.
\end{cor}
\begin{proof}
This follows from the theorem and Proposition \ref{QOrb}.
\end{proof}

\end{subsection} \begin{subsection}{Dedekind-$\eta$ Products}

In this section we review a construction that produces scalar-valued modular forms over $\widetilde{\Gamma_0(N)}$. The Dedekind-$\eta$ function is given by
\begin{eqnarray*}
\eta(\tau) = \q^{1/24}\prod\limits_{k=1}^\infty(1-\q^k)
\end{eqnarray*}
and is a weight $\frac{1}{2}$ modular form on $\MSL(\Z)$. It satisfies
\begin{eqnarray*}
\eta(\tau+1) = \expo(1/12)\eta(\tau), & &
\eta(-1/\tau) = \sqrt{-i\tau}\ \!\eta(\tau).
\end{eqnarray*}
Let $\eta_m(\tau) = \eta(m\tau)$.
\begin{theorem}[Theorem 6.2 of \cite{BorRef}]\label{EtaProd} Given the following
\begin{enumerate}
\item[1)] a lattice $L$ with the level of $\Lambda_L$ equal to $N$,
\item[2)] $r_\delta$ for $\delta \mid N$ such that $|\Lambda_L|/\prod_{\delta\mid N}\delta^{r_\delta}$ is a rational square,
\item[3)] $(1/24)\sum_{\delta\mid N}r_\delta\delta \in \Z$, and
\item[4)] $(N/24)\sum_{\delta\mid N}r_\delta/\delta \in \Z$,
\end{enumerate}
then $\prod\limits_{\delta\mid N}\eta_\delta^{r_\delta}$ is a modular form for $\widetilde{\Gamma_0(N)}$ of weight $k=\sum_\delta r_\delta/2$ and of character $\chi_{|\Lambda_L|}$ if $4 \nmid N$ and $\chi_\theta^{2k+\left(\frac{-1}{|\Lambda_L|}\right)-1}\chi_{2^{2k}|\Lambda_L|}$ if $4 \mid N$. \end{theorem}

\end{subsection} \end{section} \begin{section}{Calculating the $\kappa_\eta(m)$} \label{Chapter6}

The application of Theorem \ref{(1,2)} requires 
computing $\kappa_{\eta_-+\lambda_-}^-(m)$. This section will cover the general techniques to complete this task using the notation and results of \cite{KRY}.
Recall from Definition \ref{JSdef} that for $\mu \in L_-^\vee/L_-$ and $\psi_\mu=char(\mu+L_-)$,
\begin{eqnarray*}
E(\tau,s;\psi_\mu,+1)  = \sum\limits_{m \in \Q}E_m(\tau,s,\mu) = \sum\limits_{m\in \Q}A_\mu(s,m,v)\q^m
\end{eqnarray*}
where $\tau = u+iv$. The Fourier coefficients have Laurant expansions
\begin{eqnarray*}
A_\mu(s,m,v) = b_\mu(m,v)s + O(s^2)
\end{eqnarray*}
at $s=0$. Thus
\begin{eqnarray*}
b_\mu(m,v) &=& \ds \left\{ A_\mu(s,m,v) \right\}_{s=0} \\
&=& \q^{-m} \ds \left\{E_m(\tau,s,\mu) \right\}_{s=0}.
\end{eqnarray*}
Let $\Delta$ denote the discriminant of $\kay \simeq U$ and $h(\kay)$ its ideal class number. Following \cite{KRY}, there is a normalization $E^*_m$ which satisfies
\begin{eqnarray*}
h(\kay) \ds \left\{ E_m(\tau,s,\mu) \right\}_{s=0} = \ds \left\{ E^*_m(\tau,s,\mu) \right\}_{s=0}
\end{eqnarray*}
and a factorization of $E^*_m(\tau,s,\mu)$ into Whittaker polynomials,
\begin{eqnarray}\label{WProd}
E^*_m(\tau,s,\mu) = v^{-\frac{1}{2}}|\Delta|^{\frac{s+1}{2}}W^*_{m,\infty}(\tau,s,\mu)\prod_p W^*_{m,p}(s,\mu).
\end{eqnarray}

\begin{lemma}[Lemmas 2.4 and 2.5 of \cite{KRY}]
Suppose $U \simeq \kay$ with discriminant $\Delta$. If $L$ is unimodular (self-dual) and $m\in Q(\mu)+\Z_p$, then
\begin{eqnarray*}
W^*_{m,p}(s,\mu)=\sum\limits_{r=0}^{\text{ord}_p(m)}\left(\frac{\Delta}{p}\right)^rX^r.
\end{eqnarray*}
Thus
\begin{eqnarray*}
W^*_{m,p}(0,\mu)=\rho_p(m)=\sum\limits_{r=0}^{\text{ord}_p(m)}\left(\frac{\Delta}{p}\right)^r.
\end{eqnarray*}
If $\rho_p(m)=0$, then 
\begin{eqnarray*}
W^{*,'}_{m,p}(0,\mu)=\frac{1}{2}\log(p)(\text{ord}_p(m)+1)\rho_p(m/p).
\end{eqnarray*}
\end{lemma}
\begin{lemma}[Proposition 2.6 of \cite{KRY}]The following values are obtained at $s=0$. \nopagebreak
\begin{enumerate} 
\item[1)] $E^*_m(\tau,0,\mu) = 0$.
\item[2)] $W^*_{m,\infty}(\tau,0,\mu) = -\gamma_\infty 2 v^{\frac{1}{2}}\q^m$,
\end{enumerate}
where $\gamma_\infty$ is a local factor that will not affect later global calculations since $\prod_{p\le\infty}\gamma_p=1$.
\end{lemma}
\noindent Note that $\rho_p(m) = \rho_p(p^{\text{ord}_p(m)})$, and $\rho_p(1)=1$. Hence $W^*_{m,p}(0,\mu) \ne 1$ for only
a finite number of primes. 

\begin{theorem}
There exists a finite prime $p'$ such that $W^*_{m,p'}(0,\mu)=0$ and hence
\begin{eqnarray}
b_\mu(m,v) = \frac{-2\sqrt{|\Delta|}\gamma_\infty}{h(\kay)}\ds \left\{W^*_{m,p'}(s,\mu) \right\}_{s=0} \prod_{p \ne p'}W^*_{m,p}(0,\mu). \label{ProdFormula}
\end{eqnarray}
\end{theorem}
\noindent Note that in this case $b_\mu(m,v)$ does not depend on $v$ and thus (\ref{ProdFormula}) is equal to the limit in (\ref{kapmindef}).
Explicit formulas for $W_{m,p}(s,\mu) = \frac{W^*_{m,p}(s,\mu)}{L_p(s+1,\chi_\Delta)}$ are given in \cite{Yang}.

\end{section} \begin{section}{Examples} \label{Chapter7}
%

\begin{subsection}{$D=6$}\label{Tar}\label{Example71}
\label{Tar6} First consider the quaternion algebra ramified at the primes $2$ and $3$. 
Let $B = \left(\frac{5,6}{\Q}\right)$. 
By proposition \ref{MaxOrd}, $B$ has a maximal order given by
\begin{eqnarray}
\Oh = \Z + \left(\frac{1+\alpha}{2}\right)\Z+ \left(\frac{\alpha+\alpha\beta}{5}\right)\Z + \left(\frac{5+\alpha+5\beta+\alpha\beta}{10}\right)\Z. 
\end{eqnarray}
Further, the image of $\Gamma^*$ in $\text{PGL}_2(\R)$ is generated by three elements,
\begin{eqnarray*}
s_2 &=& -\frac{6}{5}\alpha+\beta+\frac{4}{5}\alpha\beta, \\
s_4 &=& 1-\frac{1}{5}\alpha + \frac{1}{2}\beta +\frac{3}{10}\alpha\beta, \\
s_6 &=& \frac{3}{2} - \frac{3}{10}\alpha + \frac{1}{5}\alpha\beta,
\end{eqnarray*}
which satisfy the group presentation
\begin{eqnarray*}
\langle s_2,s_4,s_6 \mid s_2^2=s_4^4=s_6^6=s_2s_4s_6=1\rangle
\end{eqnarray*}
(See Section 3.1 of \cite{Elk}). As mentioned previously, $\Xs_6$ has genus $0$ and so
there exists a parameterization $t_6: \Xs_6 \stackrel{\sim}{\rightarrow} \Pro^1$ over $\Q$. Such a map giving the isomorphism is only well-defined up to a $\text{PGL}_2$ action on $\Pro^1$. However, the map is uniquely determined once the value at three points of $\Xs_6$ are chosen. Since there are three distinguished elements of $\Gamma^*$, namely $s_2$, $s_4$, $s_6$, it is only natural to fix the value of the isomorphism at their three fixed points, $P_2$, $P_4$, $P_6$. Thus, define the map $t_6: \Xs_6 \stackrel{\sim}{\rightarrow} \Pro^1$ such that it takes on the values $0$, $1$, $\infty$ at the points $P_4$, $P_2$, $P_6$, respectively. (Warning: In \cite{Elk}, the author chooses $t_6$ to have the values $0$, $1$, $\infty$ at the points $P_2$, $P_4$, $P_6$.) This defining criteria can be expressed as
\begin{eqnarray}
\begin{array}{c}\Div(t_6) = P_4 - P_6, \\
t_6(P_2)=1.\end{array} \label{t6def}
\end{eqnarray}
Let $s^0_i$ denote the trace-$0$ part of $s_i$. Since the action of $B^\times$ factors through $\text{PGL}_2(\R)$, the fixed point of $s_i$ is the fixed point of all of $\kay_i^\times \subset B^\times$ where $\kay_i = \Q(s_i)= \Q(s^0_i)$. For the $s_i$ as above,
\begin{eqnarray}
\begin{array}{lcccr}
\kay_2 \simeq \Q(\sqrt{-6}), & & \kay_4 \simeq \Q(\sqrt{-1}), & & \kay_6 \simeq \Q(\sqrt{-3}). \label{6fields}
\end{array}
\end{eqnarray}

\begin{lemma}
The following equalities hold.
\begin{enumerate}
\item[1)] $Z(1,0;\Gamma^*) = \frac{1}{4} P_4$.
\item[2)] $Z(3,0;\Gamma^*) = \frac{1}{6} P_6$.
\end{enumerate}
\end{lemma}
\begin{proof} These identities follow from 
\begin{eqnarray*}
|\Gamma^* \bs L(1)| = |\Gamma^* \bs L(3)| = 1
\end{eqnarray*}
by Corollary \ref{A2} and that $|\text{Stab}_{\Gamma^*}(s_4^0)|=4$ and $|\text{Stab}_{\Gamma^*}(s_6^0)|=6$.
\end{proof}

\begin{prop}
\begin{eqnarray*}
\Div(t_6) &=& 4 Z(1,0;\Gamma^*) - 6 Z(3,0;\Gamma^*).
\end{eqnarray*}
\end{prop}
\noindent Hence, to use Theorem \ref{DivForm}, the input vector-valued form must have, for $m<0$,
\begin{eqnarray*}
c_0(m) = \left\{\begin{array}{ll} 2 & m=-1 \\ -3 & m=-3 \\ 0 & \text{otherwise} \end{array} \right. .
\end{eqnarray*}

\subsubsection{The Input Form}\label{Input6}
By Corollary \ref{DualSize}, $|L^\vee/L|=72$ and $N=12$. To vectorize properly, we need a form of weight $\frac{1}{2}$ and character $\chi_\theta\chi_{144}$. 

\begin{prop}Let $A_1,A_2,A_3,A_4,A_5 \in \Z$, and set
\begin{eqnarray}
r_1 &=& A_5, \label{r1}\\
r_2 &=& 16 - 12 A_1 + 36 A_2 - 9 A_3 - 14 A_4 - 6 A_5,\\
r_3 &=& -30 + 24 A_1 - 48 A_2 + 16 A_3 + 24 A_4 + 5 A_5,\\
r_4 &=& -17 + 12 A_1 - 36 A_2 + 9 A_3 + 16 A_4 + 5 A_5,\\
r_6 &=& 43 - 36 A_1 + 60 A_2 - 21 A_3 - 34 A_4 - 6 A_5,\\
r_{12} &=& -11 + 12 A_1 - 12 A_2 + 5 A_3 + 8 A_4 + A_5.
\end{eqnarray}
Then 
\begin{eqnarray} \label{r12}
\prod\limits_{\delta\mid 12}\eta_\delta^{r_\delta}
\end{eqnarray}
is a modular form for $\widetilde{\Gamma_0(12)}$ of weight $\frac{1}{2}$ and of character $\chi_\theta\chi_{144}$.
\end{prop}
\begin{proof}
One can check that the following hold.
\begin{eqnarray*}
72/\prod_{\delta\mid 12}\delta^{r_\delta} &=& (2^{A_3} 3^{A_4})^2, \\
(1/24)\sum_{\delta\mid 12}r_\delta\delta &=& A_1, \\
(1/2)\sum_{\delta\mid 12}r_\delta/\delta &=& A_2, \\
\sum_\delta r_\delta/2 &=& \frac{1}{2}.
\end{eqnarray*}
Hence, by the Theorem \ref{EtaProd}, (\ref{r12}) is a modular form for $\widetilde{\Gamma_0(12)}$ of weight $\frac{1}{2}$ and of character $\chi_\theta\chi_{144}$.
\end{proof}

Now examine the structure of such a form at the various cusps of $\widetilde{\Gamma_0(12)}$. 
Table \ref{ZeroOrd2} gives the orders of the zeroes for a form defined by (\ref{r1}-\ref{r12}), where a negative value represents a pole.
\begin{table}
\caption{Order of the zero of a form defined by (\ref{r1}-\ref{r12}) at the cusps of $\widetilde{\Gamma_0(12)}$}\label{ZeroOrd2}
\begin{center}
\begin{tabular}{c|c}
\hline
Cusp & Zero Order \\
\hline\hline
$1=0$ & $A_2/12$ \\
$1/2$ & $(15 - 12 A_1 + 28 A_2 - 8 A_3 - 12 A_4 - 4 A_5)/12$ \\
$1/3$ & $(-5 + 4 A_1 - 9 A_2 + 3 A_3 + 4 A_4 + A_5)/4$ \\
$1/4$ & $(-4 + 3 A_1 - 8 A_2 + 2 A_3 + 4 A_4 + A_5)/3$ \\
$1/6$ & $(25 -20 A_1 + 36 A_2 - 12 A_3 - 20 A_4 - 4 A_5)/4$ \\
$1/12=\infty$ & $A_1$ \\
\hline
\end{tabular}\end{center}\end{table}
To construct a form defined by (\ref{r1}-\ref{r12}) such that it has neither a pole nor a zero at $\infty$ and no pole at any finite cusp, one simply solves the following system  of inequalities over $\Z$.
\begin{eqnarray}
0 &\le& A_2, \label{InEq1}\\*
0 &\le& 15 - 12 A_1 + 28 A_2 - 8 A_3 - 12 A_4 - 4 A_5, \\*
0 &\le& -5 + 4 A_1 - 9 A_2 + 3 A_3 + 4 A_4 + A_5, \\*
0 &\le& -4 + 3 A_1 - 8 A_2 + 2 A_3 + 4 A_4 + A_5, \\*
0 &\le& 25 -20 A_1 + 36 A_2 - 12 A_3 - 20 A_4 - 4 A_5, \label{InEq2}\\*
0 &=& A_1, \label{Eq}
\end{eqnarray}
Doing so yields a unique solution
\begin{eqnarray*}
(A_1,A_2,A_3,A_4,A_5) = (0,0,1,1,-2)
\end{eqnarray*}
which produces
\begin{eqnarray*}
\psi_0=\frac{\eta_2^5}{\eta_1^2 \eta_4^2} = \sum\limits_{n \in \Z}\q^{n^2} = \theta(\tau).
\end{eqnarray*} 

Similarly a form defined by (\ref{r1}-\ref{r12}) that has a pole of order $k$ at $\infty$, but no pole at any other cusp can be found by solving the inequalities (\ref{InEq1}-\ref{InEq2}) with $A_1 = -k$ over $\Z$. For a simple pole at $\infty$, there are five such Dedekind-$\eta$ products. They are
\begin{eqnarray*}
\psi_1=\frac{\eta_2^{12} \eta_3}{\eta_1^5 \eta_4^4 \eta_6 \eta_{12}^2} &=& \frac{1}{\q}+5+O[\q],\\
\frac{\eta_2^3 \eta_4^2 \eta_6^2}{\eta_1^2 \eta_{12}^4} &=& \frac{1}{\q}+2+O[\q],\\
\frac{\eta_2^2 \eta_6^9}{\eta_1 \eta_3^3 \eta_{12}^6} &=& \frac{1}{\q}+1+O[\q],\\
\frac{\eta_2^5 \eta_3^3}{\eta_1^3 \eta_4 \eta_{12}^3} &=& \frac{1}{\q}+3+O[\q],\\
\frac{\eta_1 \eta_2^3 \eta_6^2}{\eta_3 \eta_4\eta_{12}^3} &=& \frac{1}{\q}-1+O[\q].
\end{eqnarray*}
For a triple pole, there are 35 such forms. One of them is
\begin{eqnarray*}
\psi_3 = \frac{\eta_2 \eta_3^2 \eta_4^4 \eta_6^4}{\eta_{12}^{10}} &=& \frac{1}{\q^3}-\frac{1}{\q}-2+O[\q].
\end{eqnarray*}
Thus the linear combination
\begin{eqnarray*}
f_6= -6\psi_3-2\psi_1-2\psi_0= -\frac{6}{\q^3}+\frac{4}{\q} + O[\q]
\end{eqnarray*}
is a vectorizable modular form over $\widetilde{\Gamma_0(12)}$ for $\widetilde{\Gamma_0(12)}$ of weight $\frac{1}{2}$ of character $\chi_\theta\chi_{144}$ with no poles at finite cusps.

\begin{theorem}\label{Main6}
There exists a nonzero constant $c_6$ such that
\begin{eqnarray*}
t_6 = c_6\Psi(F_{f_6})^2.
\end{eqnarray*}
\end{theorem}
\begin{proof}There is an equality of divisors
\begin{eqnarray*}
\Div(t_6) = 4 Z(1,0;\Gamma^*) - 6 Z(3,0;\Gamma^*) = \Div(\Psi(F_{f_6})^2).
\end{eqnarray*}
\end{proof}


\begin{subsubsection}{$\Delta=-24$}

In this section we calculate $\Psi(F_{f_6})(P_2)$. The result of the calculation gives the value of $c_6^{-1}$ in Theorem \ref{Main6} since by definition $t_6(P_2)=1$. Note that by (\ref{6fields}), $P_2 = \Po_{-24}$ the CM point with discriminant $-24$ on the Shimura curve $\Xs_6$.

Set $m = 1$ so that 
\begin{eqnarray*}
L &=& \Z\ell_1 + \Z\ell_2 + \Z\ell_3
\end{eqnarray*}
where
\begin{eqnarray*}
\begin{array}{lcccr}
\ell_1 = \alpha, & &
\ell_2 = \frac{\alpha+\alpha\beta}{5}, & &
\ell_3 = \frac{\beta+\alpha\beta}{2}. \end{array}
\end{eqnarray*}
Take $z=\ell_3$ so that $Q(z) = 6$. Then the negative plane is spanned by
\begin{eqnarray*}
u_1 = 2 \ell_2-\ell_3, & &
u_2 = 2 \ell_1-4 \ell_2 + 2 \ell_3,
\end{eqnarray*}
and
\begin{eqnarray*}
Q(X u_1+Yu_2) = -2(X^2+6Y^2).
\end{eqnarray*}
A basis of $L_-$ is given by
\begin{eqnarray*}
\ell_1^- = 2\ell_2-\ell_3, & &
\ell_2^- = \ell_1.
\end{eqnarray*}
The group $L/(L_-+L_+)$ has order $2$ and $\lambda=\ell_2 + (L_-+L_+)$ represents its nontrivial member. This has the decomposition
\begin{eqnarray*}
\lambda_+ = \frac{1}{2}z + L_+,& &
\lambda_- = \frac{1}{2}\ell_1^- + L_-.
\end{eqnarray*}
By Theorem \ref{(1,2)},
\begin{eqnarray}
\sum\limits_{z\in Z_{\Gamma^*}(\Q(\sqrt{-6}))}\log||\Psi(z,F_{f_6})||^2 &=& \left(\frac{-1}{4}\right)(-6\kappa_0(3)+4\kappa_0(1)). \label{key6}
\end{eqnarray}
Considering (\ref{kappasum}),
\begin{eqnarray}
\kappa_0(1) &=& \kappa_{0}^-(1),\\
\kappa_0(3) &=& \kappa_{0}^-(3) + \kappa_{\lambda_-}^-(3/2)+ \kappa_{\lambda_-}^-(3/2). \label{k03}
\end{eqnarray}
The term $\kappa_{\lambda_-}^-(3/2)$ appears twice in (\ref{k03}) due to the two values $x=\pm z/2 \in \lambda_+ + L_+ = (\frac{1}{2}+\Z)z$ that satisfy $3-Q(x) \ge 0$. 

The calculations via \cite{Yang} and Section \ref{Chapter6} yield
\begin{eqnarray*}
\kappa_0(1) &=& -6\log(2),\\
\kappa_0(3) &=& -8\log(2) -4\log(3).
\end{eqnarray*}
Thus
\begin{eqnarray*}
\sum\limits_{z\in Z_{\Gamma^*}(\Q(\sqrt{-6}))}\log||\Psi(z,F_{f_6})||^2 = -6\log(3)-6\log(2).
\end{eqnarray*}
\begin{cor} $||t_6|| = 6^{6}||\Psi(F_{f_6})^2||.$
\end{cor}
\noindent Note that we have only determined the value of $c_6$ in Theorem \ref{Main6} up to sign. This can be resolved by repeating the above computations with a Borchards form corresponding to the function $t_6-1$.
\end{subsubsection}
\begin{subsubsection}{$\Delta=-163$}

We are now able to compute the coordinates of the other rational CM points listed in Table \ref{Rat6}. We illustrate the calculations with the example of $\Delta=-163$. 

Take $z=\ell_1+14\ell_2$ so that $Q(z) = 163$. Then the negative plane is spanned by
\begin{eqnarray*}
u_1 = 42\ell_2-13\ell_3, & &
u_2 = 166 \ell_1-284 \ell_2 + 163 \ell_3,
\end{eqnarray*}
and
\begin{eqnarray*}
Q(X u_1+Yu_2) = -498(X^2+163Y^2).
\end{eqnarray*}
A basis of $L_-$ is given by
\begin{eqnarray*}
\ell_1^- = 42\ell_2-13\ell_3, \\
\ell_2^- = \ell_1-5\ell_2+2\ell_3,
\end{eqnarray*}

The group $L/(L_-+L_+)$ is cyclic of order $163$ and $\lambda=\ell_3 + (L_-+L_+)$ represents a generator. This has the decomposition
\begin{eqnarray*}
\lambda_+ &=& \frac{42}{163}z + L_+,\\
\lambda_- &=& -\frac{19}{163}\ell_1^- - \frac{42}{163} \ell_2^- + L_-.
\end{eqnarray*}
Then computations of Whittaker polynomials as before yield
\begin{eqnarray*}
\kappa_0(1) &=& -4\log(2)-11\log(3)-4\log(7)-4\log(19)-4\log(23),\\
\kappa_0(3) &=& -\frac{40}{3}\log(2)-4\log(3)-4\log(5)-4\log(11)-4\log(17).
\end{eqnarray*}
(Due to the sheer number of Whittaker polynomials required, the calculations were implemented in Mathematica.)
Thus by Theorem \ref{(1,2)} the CM point $\Po_{-163}$ with discriminant $-163$ has
\begin{eqnarray*}
||t_6(\Po_{-163})|| = \frac{3^{11}7^4 19^4 23^4}{2^{10}5^6 11^6 17^6}.
\end{eqnarray*}
Note that this proves the conjectural value given in Table 2 of \cite{Elk}. In fact, all of the conjectural values can now be algebraically proven and are given in Table \ref{Rat6}.

\end{subsubsection} \end{subsection} \begin{subsection}{$D=10$}\label{Tar10} 

\begin{subsubsection}{The Input Form}
Now consider the quaternion algebra ramified at the primes $2$ and $5$, 
$B=\left(\frac{13,10}{\Q}\right)$.
It contains the maximal order
\begin{eqnarray*}
\Oh = \Z + \left(\frac{1+\alpha}{2}\right)\Z+ \left(\frac{6\alpha+\alpha\beta}{13}\right)\Z + \left(\frac{78+6\alpha+13\beta+\alpha\beta}{26}\right)\Z.
\end{eqnarray*}
Then by Section 4.1 of \cite{Elk}, the image of $\Gamma^* \subset \text{PGL}_2(\R)$ is presented as
\begin{eqnarray*}
\langle s_2,s_2',s_2'',s_3 \mid s_2^2=s_2'^2=s_2'^2=s_3^3=s_2s_2's_2''s_3=1\rangle,
\end{eqnarray*}
with
\begin{eqnarray*}
\begin{array}{rclrcl}
s_2 &=&  -\frac{8}{13}\alpha-\frac{3}{13}\alpha\beta, &
s_2' &=& -\frac{20}{13}\alpha - \frac{1}{2}\beta - \frac{15}{26}\alpha\beta, \\
s_2'' &=& -\frac{35}{13}\alpha - \frac{1}{2}\beta - \frac{23}{26}\alpha\beta, &
s_3 &=& -\frac{1}{2} - \frac{31}{26}\alpha - \frac{5}{13}\alpha\beta.\end{array}
\end{eqnarray*}
and $\Xs_{10}$ has genus $0$. Hence, there is a map $t_{10}: \Xs_{10} \stackrel{\sim}{\rightarrow} \Pro^1$ such that
\begin{eqnarray}
\begin{array}{c}
\Div(t_{10}) = P_3 - P_2, \\
t_{10}(P_2'') = 2,
\end{array} \label{t10def}
\end{eqnarray}
where $P_2$, $P_2''$, $P_3$ are the fixed points of $s_2$, $s_2''$, $s_3$, respectively.
Again the fixed point of $s_i$ is the fixed point of all of $\kay_i^\times \subset B^\times$ where $\kay_i = \Q(s_i^0)$. Now
\begin{eqnarray*}
\begin{array}{lcccccr}
\kay_2 \simeq \Q(\sqrt{-2}), & & \kay_2' \simeq \Q(\sqrt{-10}), & & \kay_2'' \simeq \Q(\sqrt{-5}), & & \kay_3 \simeq \Q(\sqrt{-3}). \label{10fields}
\end{array}
\end{eqnarray*}

\begin{lemma}
The following equalities hold.
\begin{enumerate}
\item[1)] $Z(2,0;\Gamma^*) =  \frac{1}{2} P_2$.
\item[2)] $Z(3,0;\Gamma^*) = \frac{1}{3} P_3$.
\end{enumerate}
\end{lemma}
\begin{prop}
The following identity for $t_{10}$ holds,
\begin{eqnarray*}
\Div(t_{10}) &=& 3 Z(3,0;\Gamma^*) - 2 Z(2,0;\Gamma^*).
\end{eqnarray*}
\end{prop}
\noindent Then the same line of reasoning as in Section \ref{Input6} applied to the case $|L^\vee/L|=200$ and $N=20$ gives the following result.
\begin{theorem}\label{Main10}Let
\begin{eqnarray*}
f_{10} &=& 3\left(\frac{\eta_4^6\eta_{10}^8}{\eta_2^3\eta_5^2\eta_{20}^8}\right) - 2\left(\frac{\eta_2^3\eta_4^2\eta_{10}^2}{\eta_1^2\eta_{20}^4}\right) - 5\left(\frac{\eta_4^2\eta_{10}^6}{\eta_2\eta_5^2\eta_{20}^4}\right) +  4\left(\frac{\eta_2^5}{\eta_1^2\eta_4^2}\right) \\
&=& \frac{3}{\q^3}-\frac{2}{\q^2} + O[\q].
\end{eqnarray*}
It is a vectorizable modular form over $\widetilde{\Gamma_0(20)}$ with no poles at finite cusps. Thus
\begin{eqnarray*}
\Div(t_{10}) = 3 Z(3,0;\Gamma^*) - 2 Z(2,0;\Gamma^*) = \Div(\Psi(F_{f_{10}})^2).
\end{eqnarray*}
So again the two functions agree up to a nonzero constant,
\begin{eqnarray*}
t_{10} = c_{10}\Psi(F_{f_{10}})^2.
\end{eqnarray*}
\end{theorem}
\end{subsubsection}
\subsection{$\Delta=-20$}

To compute the constant $c_{10}$, we now consider the case of $\Delta=-20$.
Recall that $t_{10}(P_2')=2$ by definition and by (\ref{10fields}), $P_2''=\Po_{-20} \in \Xs_{10}$. Then
\begin{eqnarray*}
L &=& \Z\ell_1 + \Z\ell_2 + \Z\ell_3
\end{eqnarray*}
where
\begin{eqnarray*}
\begin{array}{lcccr}
\ell_1 = \alpha, & &
\ell_2 = \frac{6\alpha+\alpha\beta}{13}, & &
\ell_3 = \frac{\beta+\alpha\beta}{2}. \end{array}
\end{eqnarray*}
Take $z=\ell_1-3\ell_2$ so that $Q(z) = 5$. Then the negative plane is spanned by
\begin{eqnarray*}
u_1 = -\ell_2, & &
u_2 = 6 \ell_1-13 \ell_2 + 2 \ell_3,
\end{eqnarray*}
and
\begin{eqnarray*}
Q(X u_1+Yu_2) = -2(X^2+5Y^2).
\end{eqnarray*}
A basis of $L_-$ is given by
\begin{eqnarray*}
\ell_1^- = -\ell_2, & &
\ell_2^- = 3\ell_1+\ell_3,
\end{eqnarray*}

In this case the quotient $L/(L_-+L_+)$ is trivial. Theorem \ref{(1,2)} yields,
\begin{eqnarray*}
\sum\limits_{z\in Z_{\Gamma^*}(\Q(\sqrt{-5}))}\log||\Psi(z,F_{f_{10}})||^2 &=& \left(\frac{-1}{4}\right)(3\kappa_0(3)-2\kappa_0(2)) = 3 \log(2).
\end{eqnarray*}
Thus
\begin{eqnarray*}
||\Psi(P_2'',F_{f_{10}})^2|| = 2^{3}.
\end{eqnarray*}
Since $t_{10}(P_2'')=2$,
\begin{eqnarray*}
||t_{10}|| = 2^{-2}||\Psi(F_{f_{10}})^2||.
\end{eqnarray*}

\subsection{$\Delta=-68$}

Again, we are now able to compute the coordinates of the other rational CM points for $\Xs_{10}$ listed in Table \ref{Rat10}. Moreover, we are also capable of calculating the norms of irrational CM points. As an example, we compute the norm of the irrational CM point with discriminant $-68$.

Take $z=7\ell_1-13\ell_2+\ell_3$ so that $Q(z) = 17$. Then the negative plane is spanned by
\begin{eqnarray*}
u_1 = -35\ell_2+11\ell_3, & &
u_2 = 534 \ell_1-1067 \ell_2 + 137 \ell_3,
\end{eqnarray*}
and
\begin{eqnarray*}
Q(X u_1+Yu_2) = -2670(X^2+17Y^2).
\end{eqnarray*}
A basis of $L_-$ is given by
\begin{eqnarray*}
\ell_1^- = -35\ell_2+11\ell_3, & &
\ell_2^- = \ell_1+2\ell_2-\ell_3,
\end{eqnarray*}

The group $L/(L_-+L_+)$ is cyclic of order $17$ and is generated by $\lambda=\ell_3 + (L_-+L_+)$. This has the decomposition
\begin{eqnarray*}
\lambda_+ = \frac{-35}{17}z + L_+, & &
\lambda_- = \frac{27}{17}\ell_1^- +\frac{245}{17} \ell_2^- + L_-.
\end{eqnarray*}
Then computations as before yield 
\begin{eqnarray*}
\kappa_0(1) &=& -6\log(2)-6\log(5),\\
\kappa_0(3) &=& -8\log(2)-\frac{14}{3}\log(5).
\end{eqnarray*}
This time the CM point with discriminant $-68$ is irrational, and Theorem \ref{(1,2)} gives its norm (after renormalization) as 
\begin{eqnarray*}
\prod_{z \in Z(-68)}||\Psi(z,F_{f_{10}})||^2 = 2^2\cdot5.
\end{eqnarray*}

\end{subsection} \end{section} \begin{section}{Tables}

\begin{subsection}{$D=6$}

\subsection{Coordinates of Rational CM Points on $\Xs_6$}
The following table gives the values of $t_6$ (as defined by (\ref{t6def})) at the rational CM points of $\Xs_6$. These values verify the values in Table 2 of \cite{Elk}. Denote $t_6(P_{CM}) = (r:s)$. 
\begin{center} 
\begin{longtable}{|c|c|c|c|}
\caption{Coordinates of Rational CM Points on $\Xs_6$}\label{Rat6} \\
\hline
$\Delta$ & $r$ & $s$ & Proved in \cite{Elk}\\
\hline\endhead
$-3$   & $1$ & $0$                  & Y  \\
$-4$   & $0$ & $1$                 & Y  \\
$-24$  & $1$ & $1$                 & Y  \\ 
$-40$  & $3^7$ & $5^3$             & Y  \\
$-52$  & $2^23^7$ & $5^6$          & Y  \\
$-19$  & $3^7$ & $2^{10}$          & Y  \\
$-84$  & $-2^27^2$ & $3^3$           & Y \\
$-88$  & $3^77^4$ & $5^611^3$       & Y \\
$-100$ & $2^43^77^45$ & $11^6$           & Y \\
$-120$ & $7^4$ & $3^35^3$           & Y \\
$-132$ & $2^411^2$ & $5^6$ & Y \\
$-148$ & $2^23^77^411^4$ & $5^617^6$ & N \\
$-168$ & $-7^211^4$ & $5^6$ & Y \\
$-43$  & $3^77^4$ & $2^{10}5^6$    & Y  \\
$-51$  & $-7^4$ & $2^{10}$          & Y  \\
$-228$ & $2^67^419^2$ & $3^65^6$ & N \\
$-232$ & $3^77^411^419^4$ & $5^623^629^3$ & N \\
$-67$  & $3^77^411^4$ & $2^{16}5^6$ & N \\
$-75$  & $11^4$ & $2^{10}3^35$ & Y \\
$-312$ & $7^423^4$ & $5^611^6$ & Y \\
$-372$ & $-2^27^419^431^2$ & $3^35^611^6$ & N \\
$-408$ & $-7^411^431^4$ & $3^65^617^3$ & N \\
$-123$ & $-7^419^4$ & $2^{10}5^6$ & N \\
$-147$ & $-11^423^4$ & $2^{10}3^35^67$ & Y \\
$-163$ & $3^{11}7^419^423^4$ & $2^{10}5^611^617^6$ & N \\
$-708$ & $2^87^411^447^459^2$ & $5^617^629^6$ & N \\
$-267$ & $-7^431^443^4$ & $2^{16}5^611^6$ & N \\
\hline
\end{longtable}\end{center}

\subsection{Norms of CM Points on $\Xs_6$ for $0< -d \le 250$}
Here we give the norms for all CM points of fundamental discriminant $\Delta = d$ or $4d$ for $0< -d \le 250$. This cut-off is arbitrary. It is also only for implementation reasons that we only compute for fundamental discriminants (i.e. $d$ squarefree). 
\begin{center} 
\begin{longtable}{|c|c|c|}
\caption{Norms of CM Points on $\Xs_6$ for $0< -d \le 250$}\label{Norm6} \\
\hline
$\Delta$ &$|t_6(\Po_\Delta)|$ & $|(1-t_6)(\Po_\Delta)|$ \\
\hline\endhead
$-40$ & $\frac{3^{7}}{5^{3}}$ & $\frac{2^{3}17^{2}}{5^{3}}$ \\
$-52$ & $\frac{2^{2}3^{7}}{5^{6}}$ & $\frac{13^{1}23^{2}}{5^{6}}$ \\
$-19$ & $\frac{3^{7}}{2^{10}}$ & $\frac{13^{2}19^{1}}{2^{10}}$ \\
$-84$ & $\frac{2^{2}7^{2}}{3^{3}}$ & $\frac{13^{2}}{3^{3}}$ \\
$-88$ & $\frac{3^{7}7^{4}}{5^{6}11^{3}}$ & $\frac{2^{5}17^{2}41^{2}}{5^{6}11^{3}}$ \\
$-120$ & $\frac{7^{4}}{3^{3}5^{3}}$ & $\frac{2^{4}19^{2}}{3^{3}5^{3}}$ \\
$-132$ & $\frac{2^{4}11^{2}}{5^{6}}$ & $\frac{3^{4}13^{2}}{5^{6}}$ \\
$-136$ & $\frac{3^{14}}{11^{2}17^{3}}$ & $\frac{2^{6}13^{4}41^{2}}{11^{6}17^{2}}$ \\
$-148$ & $\frac{2^{2}3^{7}7^{4}11^{4}}{5^{6}17^{6}}$ & $\frac{13^{2}37^{1}47^{2}71^{2}}{5^{6}17^{6}}$ \\
$-168$ & $\frac{7^{2}11^{4}}{5^{6}}$ & $\frac{2^{3}3^{5}19^{2}}{5^{6}}$ \\
$-43$ & $\frac{3^{7}7^{4}}{2^{10}5^{6}}$ & $\frac{19^{2}37^{2}43^{1}}{2^{10}5^{6}}$ \\
$-184$ & $\frac{3^{14}7^{8}}{17^{6}23^{3}}$ & $\frac{2^{8}13^{4}89^{2}}{17^{4}23^{2}}$ \\
$-51$ & $\frac{7^{4}}{2^{10}}$ & $\frac{3^{4}17^{1}}{2^{10}}$ \\
$-228$ & $\frac{2^{6}7^{4}19^{2}}{3^{6}5^{6}}$ & $\frac{13^{2}17^{2}37^{2}}{3^{6}5^{6}}$ \\
$-232$ & $\frac{3^{7}7^{4}11^{4}19^{4}}{5^{6}23^{6}29^{3}}$ & $\frac{2^{3}13^{2}17^{2}41^{2}89^{2}113^{2}}{5^{6}23^{6}29^{3}}$ \\
$-244$ & $\frac{2^{6}3^{21}19^{4}}{17^{6}29^{6}}$ & $\frac{19^{4}37^{2}47^{2}61^{1}}{17^{2}29^{6}}$ \\
$-264$ & $\frac{19^{4}}{3^{9}11^{1}}$ & $\frac{2^{6}19^{2}43^{2}}{3^{9}11^{3}}$ \\
$-67$ & $\frac{3^{7}7^{4}11^{4}}{2^{16}5^{6}}$ & $\frac{13^{2}43^{2}61^{2}67^{1}}{2^{16}5^{6}}$ \\
$-276$ & $\frac{2^{4}23^{2}}{11^{2}}$ & $\frac{3^{8}23^{1}37^{2}}{11^{6}}$ \\
$-280$ & $\frac{3^{14}7^{4}23^{4}}{5^{6}11^{2}29^{6}}$ & $\frac{2^{12}13^{4}23^{2}113^{2}137^{2}}{5^{6}11^{6}29^{6}}$ \\
$-292$ & $\frac{2^{10}3^{14}19^{4}}{5^{12}23^{2}}$ & $\frac{13^{4}17^{4}19^{2}67^{2}71^{2}}{5^{12}23^{6}}$ \\
$-312$ & $\frac{7^{4}23^{4}}{5^{6}11^{6}}$ & $\frac{2^{4}3^{5}13^{1}17^{2}43^{2}}{5^{6}11^{6}}$ \\
$-328$ & $\frac{3^{18}11^{8}19^{4}}{5^{12}17^{6}41^{3}}$ & $\frac{2^{6}19^{2}23^{4}89^{2}137^{2}}{5^{12}17^{4}41^{2}}$ \\
$-340$ & $\frac{2^{4}3^{18}7^{8}23^{4}}{5^{6}29^{6}41^{6}}$ & $\frac{13^{4}17^{2}19^{4}23^{2}61^{2}167^{2}}{5^{6}29^{6}41^{6}}$ \\
$-91$ & $\frac{3^{14}7^{4}}{2^{26}11^{2}}$ & $\frac{13^{2}17^{4}37^{2}67^{2}}{2^{26}11^{6}}$ \\
$-372$ & $\frac{2^{2}7^{4}19^{4}31^{2}}{3^{3}5^{6}11^{6}}$ & $\frac{13^{2}23^{2}37^{2}61^{2}}{3^{3}5^{6}11^{6}}$ \\
$-376$ & $\frac{3^{28}31^{4}}{23^{2}41^{6}47^{3}}$ & $\frac{2^{16}37^{4}113^{2}}{23^{2}41^{4}47^{2}}$ \\
$-388$ & $\frac{2^{14}3^{18}31^{4}}{5^{12}11^{2}47^{6}}$ & $\frac{13^{4}17^{4}43^{2}167^{2}191^{2}}{5^{12}11^{6}47^{4}}$ \\
$-408$ & $\frac{7^{4}11^{4}31^{4}}{3^{6}5^{6}17^{3}}$ & $\frac{2^{6}13^{2}19^{2}43^{2}67^{2}}{3^{6}5^{6}17^{3}}$ \\
$-420$ & $\frac{2^{12}7^{4}23^{4}}{5^{6}17^{6}}$ & $\frac{3^{8}23^{2}61^{2}}{5^{6}17^{4}}$ \\
$-424$ & $\frac{3^{25}7^{12}19^{4}}{29^{6}47^{6}53^{3}}$ & $\frac{2^{9}13^{6}19^{4}37^{4}41^{2}137^{2}}{29^{6}47^{6}53^{3}}$ \\
$-436$ & $\frac{2^{6}3^{21}7^{12}31^{4}}{17^{6}41^{6}53^{6}}$ & $\frac{13^{6}43^{4}71^{2}109^{1}191^{2}}{17^{2}41^{6}53^{6}}$ \\
$-456$ & $\frac{7^{8}19^{2}}{11^{2}17^{6}}$ & $\frac{2^{6}3^{9}19^{1}67^{2}}{11^{6}17^{4}}$ \\
$-115$ & $\frac{3^{14}19^{4}}{2^{20}5^{6}11^{2}}$ & $\frac{13^{4}19^{2}23^{2}61^{2}109^{2}}{2^{20}5^{6}11^{6}}$ \\
$-472$ & $\frac{3^{21}19^{4}23^{4}31^{4}}{5^{18}53^{6}59^{3}}$ & $\frac{2^{17}19^{4}23^{4}47^{2}89^{2}233^{2}}{5^{18}53^{6}59^{3}}$ \\
$-123$ & $\frac{7^{4}19^{4}}{2^{10}5^{6}}$ & $\frac{3^{4}13^{2}23^{2}41^{1}}{2^{10}5^{6}}$ \\
$-516$ & $\frac{2^{14}31^{4}43^{2}}{3^{12}17^{6}}$ & $\frac{37^{2}41^{2}61^{2}}{3^{12}17^{2}}$ \\
$-520$ & $\frac{3^{18}7^{8}19^{4}43^{4}}{5^{6}11^{2}41^{6}59^{6}}$ & $\frac{2^{6}13^{2}17^{4}19^{2}37^{4}113^{2}233^{2}257^{2}}{5^{6}11^{6}41^{6}59^{6}}$\\
$-532$ & $\frac{2^{4}3^{14}7^{4}11^{8}23^{4}43^{4}}{5^{12}29^{6}53^{6}}$ & $\frac{17^{4}19^{2}23^{2}37^{2}109^{2}191^{2}239^{2}263^{2}}{5^{12}29^{6}53^{6}}$ \\ 
$-552$ & $\frac{19^{4}43^{4}}{3^{6}5^{12}23^{1}}$ & $\frac{2^{6}13^{4}19^{2}43^{2}67^{2}}{3^{6}5^{12}23^{3}}$ \\
$-139$ & $\frac{3^{21}19^{4}23^{4}}{2^{36}17^{6}}$ & $\frac{19^{4}23^{4}43^{2}139^{1}}{2^{36}17^{2}}$ \\
$-564$ & $\frac{2^{4}7^{8}47^{2}}{11^{2}23^{6}}$ & $\frac{3^{8}13^{4}17^{4}47^{1}}{11^{6}23^{4}}$ \\
$-568$ & $\frac{3^{14}7^{8}23^{4}31^{4}}{5^{12}17^{6}47^{2}71^{3}}$ & $\frac{2^{8}19^{4}23^{2}41^{2}137^{2}257^{2}281^{2}}{5^{12}17^{4}47^{6}71^{2}}$\\
$-580$ & $\frac{2^{20}3^{32}43^{4}47^{4}}{5^{12}59^{6}71^{6}}$ & $\frac{13^{8}41^{4}43^{4}47^{2}139^{2}263^{2}}{5^{12}59^{6}71^{6}}$ \\
$-616$ & $\frac{3^{28}7^{8}43^{4}}{11^{2}23^{2}53^{6}71^{6}}$ & $\frac{2^{12}13^{8}37^{2}61^{4}233^{2}281^{2}}{11^{6}23^{2}53^{6}71^{6}}$ \\
$-628$ & $\frac{2^{6}3^{21}19^{4}31^{4}47^{4}}{5^{18}11^{8}41^{6}}$ & $\frac{19^{4}61^{2}71^{2}157^{1}167^{2}239^{2}311^{2}}{5^{18}11^{12}41^{2}}$ \\
$-163$ & $\frac{3^{11}7^{4}19^{4}23^{4}}{2^{10}5^{6}11^{6}17^{6}}$ & $\frac{13^{2}67^{2}109^{2}139^{2}157^{2}163^{1}}{2^{10}5^{6}11^{6}17^{6}}$ \\
$-660$ & $\frac{2^{4}7^{8}11^{4}43^{4}}{3^{12}5^{6}23^{6}}$ & $\frac{17^{4}47^{2}61^{2}109^{2}}{3^{12}5^{6}23^{4}}$ \\
$-664$ & $\frac{3^{39}47^{4}}{29^{6}59^{6}83^{3}}$ & $\frac{2^{27}37^{4}47^{4}61^{4}71^{2}89^{2}257^{2}}{11^{12}29^{6}59^{6}83^{3}}$\\
$-696$ & $\frac{11^{2}23^{4}31^{4}}{3^{15}29^{3}}$ & $\frac{2^{12}13^{6}23^{4}67^{2}}{3^{15}11^{6}29^{3}}$\\
$-708$ & $\frac{2^{8}7^{4}11^{4}47^{4}59^{2}}{5^{6}17^{6}29^{6}}$ & $\frac{3^{4}13^{2}19^{2}23^{2}37^{2}41^{2}109^{2}}{5^{6}17^{6}29^{6}}$ \\
$-712$ & $\frac{3^{32}19^{8}43^{4}59^{4}}{5^{24}83^{6}89^{3}}$ & $\frac{2^{12}19^{4}41^{2}43^{2}47^{4}113^{2}281^{2}353^{2}}{5^{24}83^{6}89^{2}}$  \\
$-724$ & $\frac{2^{10}3^{39}59^{4}}{17^{6}53^{6}89^{6}}$ & $\frac{17^{2}41^{4}43^{4}67^{4}157^{2}181^{1}359^{2}}{11^{12}53^{6}89^{6}}$ \\
$-744$ & $\frac{7^{12}31^{2}59^{4}}{23^{6}29^{6}}$ & $\frac{2^{9}3^{13}41^{2}43^{2}}{23^{2}29^{6}}$ \\
$-187$ & $\frac{3^{18}11^{4}31^{4}}{2^{20}5^{12}23^{2}}$ & $\frac{13^{4}17^{2}19^{4}37^{2}163^{2}181^{2}}{2^{20}5^{12}23^{6}}$ \\
$-760$ & $\frac{3^{14}7^{8}11^{8}23^{4}31^{4}47^{4}}{5^{6}41^{6}71^{6}89^{6}}$ & $\frac{2^{8}13^{4}17^{4}19^{2}23^{2}47^{2}61^{4}137^{2}233^{2}353^{2}}{5^{6}41^{6}71^{6}89^{6}}$\\
$-772$ & $\frac{2^{14}3^{14}7^{8}31^{4}43^{4}}{5^{12}23^{6}59^{2}83^{6}}$ & $\frac{13^{4}17^{4}43^{2}139^{2}239^{2}311^{2}359^{2}383^{2}}{5^{12}23^{4}59^{6}83^{6}}$\\
$-195$ & $\frac{19^{4}31^{4}}{2^{26}5^{6}}$ & $\frac{3^{8}13^{2}19^{2}47^{2}}{2^{26}5^{6}}$\\
$-804$ & $\frac{2^{18}11^{2}19^{4}67^{2}}{3^{12}29^{6}}$ & $\frac{13^{6}17^{6}19^{4}109^{2}}{3^{12}11^{6}29^{6}}$ \\
$-808$ & $\frac{3^{25}23^{4}31^{4}59^{4}67^{4}}{5^{18}11^{8}47^{6}101^{3}}$ & $\frac{2^{9}13^{6}23^{4}37^{2}41^{2}89^{2}257^{2}401^{2}}{5^{18}11^{12}47^{2}101^{3}}$\\
$-820$ & $\frac{2^{8}3^{28}7^{16}47^{4}67^{4}}{5^{12}29^{6}89^{6}101^{6}}$ & $\frac{37^{4}41^{2}47^{2}67^{4}109^{2}167^{2}181^{2}263^{2}383^{2}}{5^{12}29^{6}89^{6}101^{6}}$\\
$-840$ & $\frac{7^{4}43^{4}67^{4}}{3^{12}5^{6}11^{2}17^{6}}$ & $\frac{2^{6}13^{4}19^{4}23^{4}139^{2}}{3^{12}5^{6}11^{6}17^{4}}$ \\
$-211$ & $\frac{3^{25}7^{12}31^{4}}{2^{36}17^{6}23^{6}}$ & $\frac{41^{4}61^{2}157^{2}211^{1}}{2^{36}17^{2}23^{2}}$ \\
$-852$ & $\frac{2^{4}59^{4}71^{2}}{5^{12}11^{2}23^{2}}$ & $\frac{3^{8}13^{4}19^{4}47^{2}61^{2}71^{1}}{5^{12}11^{6}23^{6}}$ \\
$-856$ & $\frac{3^{21}7^{12}11^{2}19^{4}31^{4}71^{4}}{53^{6}83^{6}101^{6}107^{3}}$ & $\frac{2^{19}13^{6}17^{6}19^{4}37^{4}281^{2}353^{2}401^{2}}{11^{6}53^{6}83^{6}101^{6}107^{3}}$\\
$-868$ & $\frac{2^{28}3^{28}7^{8}67^{4}}{5^{24}71^{2}107^{6}}$ & $\frac{37^{4}41^{4}67^{2}163^{2}191^{2}211^{2}359^{2}431^{2}}{5^{24}71^{4}107^{6}}$\\
$-219$ & $\frac{7^{8}23^{4}}{2^{26}3^{9}}$ & $\frac{13^{4}23^{2}41^{2}71^{2}}{2^{26}3^{9}}$ \\
$-888$ & $\frac{31^{4}47^{4}71^{4}}{5^{18}29^{6}}$ & $\frac{2^{12}3^{13}37^{1}41^{2}67^{2}139^{2}}{5^{18}29^{6}}$ \\
$-904$ & $\frac{3^{32}7^{16}19^{4}67^{4}}{17^{6}59^{2}89^{6}107^{6}113^{3}}$ & $\frac{2^{12}13^{8}19^{6}43^{2}61^{4}449^{2}}{59^{6}89^{4}107^{6}113^{2}}$ \\
$-916$ & $\frac{2^{10}3^{35}19^{4}43^{4}71^{4}}{41^{6}101^{6}113^{6}}$ & $\frac{13^{10}19^{8}43^{4}229^{1}311^{2}383^{2}431^{2}}{11^{12}41^{2}101^{6}113^{6}}$\\
$-235$ & $\frac{3^{14}7^{8}19^{4}31^{4}}{2^{20}5^{6}11^{2}29^{6}}$ & $\frac{17^{4}19^{2}47^{2}139^{2}181^{2}211^{2}229^{2}}{2^{20}5^{6}11^{6}29^{6}}$  \\
$-948$ & $\frac{2^{6}19^{4}31^{4}67^{4}79^{2}}{3^{15}5^{18}}$ & $\frac{19^{4}37^{2}47^{2}71^{2}109^{2}157^{2}}{3^{15}5^{18}}$ \\
$-952$ & $\frac{3^{28}7^{8}23^{2}71^{4}79^{4}}{5^{24}17^{6}113^{6}}$ & $\frac{2^{16}43^{4}47^{4}71^{2}233^{2}401^{2}449^{2}}{5^{24}17^{4}23^{4}113^{4}}$  \\
$-964$ & $\frac{2^{34}3^{42}59^{4}79^{4}}{17^{12}47^{2}83^{6}107^{6}}$ & $\frac{13^{12}37^{4}67^{4}239^{2}479^{2}}{17^{4}47^{2}83^{6}107^{6}}$ \\
$-984$ & $\frac{7^{12}11^{2}79^{4}}{3^{12}23^{6}41^{3}}$ & $\frac{2^{16}37^{2}43^{2}139^{2}163^{2}}{3^{12}11^{6}23^{2}41^{3}}$ \\
$-996$ & $\frac{2^{16}7^{12}71^{4}83^{2}}{17^{6}29^{6}41^{6}}$ & $\frac{3^{14}13^{6}47^{2}157^{2}}{17^{2}29^{6}41^{6}}$ \\
\hline
\end{longtable}\end{center}

\end{subsection} \begin{subsection}{$D=10$}
\subsection{Coordinates of Rational CM Points on $\Xs_{10}$}
The following table gives the values of $t_{10}$ (as defined by (\ref{t10def})) at the rational CM points of $\Xs_{10}$. These values verify the values in Table 4 of \cite{Elk}. Again denote $t_{10}(P_{CM}) = (r:s)$. 
\begin{center} 
\begin{longtable}{|c|c|c|c|c|}
\caption{Coordinates of Rational CM Points on $\Xs_{10}$}\label{Rat10} \\
\hline
$\Delta$ & $r$ & $s$ & Proved in \cite{Elk} \\
\hline\endhead
$-3$   & $0$ & $1$                            & Y  \\
$-8$   & $1$ & $0$                            & Y  \\
$-20$  & $2$ & $1$                            & Y  \\ 
$-40$  & $3^3$ & $1$                          & Y  \\
$-52$  & $-2^{1}3^3$ & $5^2$                  & N  \\
$-72$  & $5^3$ & $3^{1}7^2$                  & Y  \\
$-120$  & $-3^3$ & $7^2$                       & Y  \\
$-88$  & $3^35^3$ & $2^{1}7^2$               & N  \\
$-27$  & $-2^63$ & $5^2$                       & Y  \\
$-35$  & $2^6$ & $7$                          & Y  \\
$-148$  & $2^{1}3^311^3$ & $5^27^213^2$      & N  \\
$-43$  & $2^63^3$ & $5^27^2$                  & N  \\
$-180$  & $-2^{1}11^3$ & $13^2$               & Y  \\
$-232$  & $3^311^317^3$ & $2^25^27^223^2$     & N  \\
$-67$  & $-2^63^35^3$ & $7^213^2$              & N  \\
$-280$  & $3^311^3$ & $2^{1}7^{1}23^2$      & N  \\
$-340$  & $2^{1}3^323^3$ & $7^229^2$         & N  \\
$-115$  & $2^93^3$ & $13^223$                 & N  \\
$-520$  & $3^329^3$ & $2^37^213^{1}47^2$     & N  \\
$-163$  & $-2^93^35^311^3$ & $7^213^229^231^2$ & N  \\
$-760$  & $3^317^347^3$ & $7^231^271^2$       & N  \\
$-235$  & $2^63^317^3$ & $7^237^247$          & N  \\
\hline
\end{longtable}\end{center}

\subsection{Norms of CM Points on $\Xs_{10}$ for $0< -d \le 250$}
\begin{center} 
\begin{longtable}{|c|c|c|}
\caption{Norms of CM Points on $\Xs_{10}$ for $0< -d \le 250$}\label{Norm10} \\
\hline
$\Delta$ &$|t_{10}(\Po_\Delta)|$ & $|(2-t_{10})(\Po_\Delta)|$ \\
\hline\endhead
$-40$ & $\frac{3^{3}}{1}$ & $\frac{5^{2}}{1}$ \\
$-52$ & $\frac{2^{1}3^{3}}{5^{2}}$ & $\frac{2^{3}13^{1}}{5^{2}}$ \\
$-68$ & $\frac{2^{2}5^{1}}{1}$ & $\frac{2^{4}17^{1}}{5^{2}}$ \\
$-88$ & $\frac{3^{3}5^{3}}{2^{1}7^{2}}$ & $\frac{11^{1}17^{2}}{2^{1}7^{2}}$ \\
$-120$ & $\frac{3^{3}}{7^{2}}$ & $\frac{5^{3}}{7^{2}}$ \\
$-132$ & $\frac{2^{2}3^{6}5^{1}}{13^{2}}$ & $\frac{2^{4}11^{2}}{5^{2}}$ \\
$-35$ & $\frac{2^{6}}{7^{1}}$ & $\frac{2^{1}5^{2}}{7^{1}}$ \\
$-148$ & $\frac{2^{1}3^{3}11^{3}}{5^{2}7^{2}13^{2}}$ & $\frac{2^{5}17^{2}37^{1}}{5^{2}7^{2}13^{2}}$ \\
$-152$ & $\frac{11^{3}}{2^{1}5^{1}}$ & $\frac{11^{4}19^{1}}{2^{1}5^{4}}$ \\
$-168$ & $\frac{3^{6}11^{3}}{2^{2}5^{4}7^{2}}$ & $\frac{11^{2}37^{2}}{2^{2}5^{4}7^{2}}$ \\
$-43$ & $\frac{2^{6}3^{3}}{5^{2}7^{2}}$ & $\frac{2^{1}19^{2}}{5^{2}7^{2}}$ \\
$-212$ & $\frac{2^{3}5^{4}11^{3}}{7^{6}}$ & $\frac{2^{11}11^{4}53^{1}}{5^{2}7^{6}}$ \\
$-228$ & $\frac{2^{2}3^{6}5^{1}17^{3}}{7^{4}13^{2}}$ & $\frac{2^{4}19^{2}37^{2}}{5^{2}7^{4}}$ \\
$-232$ & $\frac{3^{3}11^{3}17^{3}}{2^{2}5^{2}7^{2}23^{2}}$ & $\frac{13^{2}19^{2}53^{2}}{2^{2}5^{2}7^{2}23^{2}}$ \\
$-248$ & $\frac{5^{2}17^{3}}{2^{2}23^{2}}$ & $\frac{17^{2}19^{4}31^{1}}{2^{2}5^{4}23^{2}}$ \\
$-260$ & $\frac{2^{2}17^{3}}{7^{4}13^{1}}$ & $\frac{2^{4}5^{4}}{7^{4}}$ \\
$-67$ & $\frac{2^{6}3^{3}5^{3}}{7^{2}13^{2}}$ & $\frac{2^{1}11^{2}31^{2}}{7^{2}13^{2}}$ \\
$-280$ & $\frac{3^{3}11^{3}}{2^{1}7^{1}23^{2}}$ & $\frac{5^{3}13^{2}}{2^{1}7^{1}23^{2}}$ \\
$-292$ & $\frac{2^{2}3^{6}5^{1}17^{3}}{13^{4}29^{2}}$ & $\frac{2^{4}17^{2}53^{2}73^{1}}{5^{2}13^{4}29^{2}}$ \\
$-308$ & $\frac{2^{4}5^{2}11^{3}23^{3}}{7^{4}29^{2}}$ & $\frac{2^{14}11^{2}19^{4}}{5^{4}7^{4}29^{2}}$ \\
$-312$ & $\frac{3^{6}17^{3}23^{3}}{2^{2}5^{4}7^{4}31^{2}}$ & $\frac{11^{4}13^{2}73^{2}}{2^{2}5^{4}7^{4}31^{2}}$ \\
$-328$ & $\frac{3^{6}5^{1}23^{1}}{2^{3}31^{2}}$ & $\frac{11^{4}37^{2}}{2^{3}5^{2}23^{2}}$ \\
$-83$ & $\frac{2^{18}}{5^{1}13^{2}}$ & $\frac{2^{3}13^{2}19^{2}}{5^{4}}$ \\
$-340$ & $\frac{2^{1}3^{3}23^{3}}{7^{2}29^{2}}$ & $\frac{2^{3}5^{2}13^{2}17^{1}}{7^{2}29^{2}}$ \\
$-372$ & $\frac{2^{2}3^{9}11^{3}23^{3}}{5^{4}7^{4}13^{2}37^{2}}$ & $\frac{2^{8}11^{2}31^{2}73^{2}}{5^{4}7^{4}37^{2}}$ \\
$-388$ & $\frac{2^{2}3^{6}17^{3}29^{1}}{5^{4}13^{2}37^{2}}$ & $\frac{2^{4}11^{4}17^{2}97^{1}}{5^{4}29^{2}37^{2}}$ \\
$-408$ & $\frac{3^{9}5^{1}11^{3}29^{3}}{7^{4}13^{4}31^{2}}$ & $\frac{11^{2}17^{2}19^{4}97^{2}}{5^{2}7^{4}13^{4}31^{2}}$ \\
$-420$ & $\frac{2^{2}3^{6}29^{3}}{7^{2}37^{2}}$ & $\frac{2^{4}5^{6}17^{2}}{7^{2}37^{2}}$ \\
$-107$ & $\frac{2^{21}}{5^{1}7^{6}}$ & $\frac{2^{3}17^{4}31^{2}}{5^{4}7^{6}}$ \\
$-440$ & $\frac{11^{3}23^{3}}{2^{2}13^{4}}$ & $\frac{5^{6}11^{1}}{2^{2}13^{2}}$ \\
$-452$ & $\frac{2^{4}11^{6}17^{3}29^{1}}{5^{3}7^{8}}$ & $\frac{2^{8}11^{4}17^{2}31^{4}113^{1}}{5^{6}7^{8}29^{2}}$ \\
$-115$ & $\frac{2^{9}3^{3}}{13^{2}23^{1}}$ & $\frac{2^{1}5^{2}11^{2}}{13^{2}23^{1}}$ \\
$-472$ & $\frac{3^{9}5^{4}29^{3}}{2^{1}23^{2}31^{2}47^{2}}$ & $\frac{19^{2}59^{1}73^{2}113^{2}}{2^{1}5^{2}23^{2}31^{2}47^{2}}$ \\
$-488$ & $\frac{11^{9}17^{3}}{2^{4}13^{4}47^{2}}$ & $\frac{11^{4}13^{2}17^{4}}{2^{4}5^{6}47^{2}}$ \\
$-123$ & $\frac{2^{15}3^{6}5^{1}}{7^{4}23^{2}}$ & $\frac{2^{2}13^{4}59^{2}}{5^{2}7^{4}23^{2}}$ \\
$-520$ & $\frac{3^{3}29^{3}}{2^{3}7^{2}13^{1}47^{2}}$ & $\frac{5^{4}11^{2}17^{2}}{2^{3}7^{2}13^{1}47^{2}}$ \\
$-532$ & $\frac{2^{2}3^{6}5^{6}11^{3}23^{1}}{7^{2}29^{2}37^{2}53^{2}}$ & $\frac{2^{10}11^{2}19^{2}113^{2}}{7^{2}23^{2}29^{2}37^{2}}$ \\
$-548$ & $\frac{2^{4}11^{6}29^{3}41^{3}}{5^{3}7^{8}13^{2}53^{2}}$ & $\frac{2^{8}11^{4}13^{4}19^{4}137^{1}}{5^{6}7^{8}53^{2}}$ \\
$-552$ & $\frac{3^{15}5^{2}41^{3}}{2^{6}13^{6}31^{2}}$ & $\frac{19^{4}31^{2}59^{2}}{2^{6}5^{4}13^{4}}$ \\
$-568$ & $\frac{3^{6}17^{3}23^{1}41^{3}}{5^{4}7^{4}47^{2}}$ & $\frac{17^{2}31^{2}71^{1}97^{2}137^{2}}{5^{4}7^{4}23^{2}47^{2}}$ \\
$-580$ & $\frac{2^{2}3^{6}41^{3}}{13^{2}29^{1}53^{2}}$ & $\frac{2^{4}5^{7}}{29^{1}53^{2}}$ \\
$-155$ & $\frac{2^{12}11^{3}}{7^{4}31^{1}}$ & $\frac{2^{2}5^{6}11^{2}}{7^{4}31^{1}}$ \\
$-628$ & $\frac{2^{3}3^{9}5^{4}11^{3}47^{3}}{29^{2}31^{2}53^{2}61^{2}}$ & $\frac{2^{13}11^{4}19^{2}137^{2}157^{1}}{5^{2}29^{2}31^{2}53^{2}61^{2}}$ \\
$-632$ & $\frac{5^{2}11^{6}41^{3}47^{1}}{2^{2}7^{8}13^{4}}$ & $\frac{11^{4}17^{6}79^{1}113^{2}}{2^{2}5^{4}7^{8}47^{2}}$ \\
$-163$ & $\frac{2^{9}3^{3}5^{3}11^{3}}{7^{2}13^{2}29^{2}31^{2}}$ & $\frac{2^{1}19^{2}59^{2}79^{2}}{7^{2}13^{2}29^{2}31^{2}}$ \\
$-660$ & $\frac{2^{2}3^{9}47^{3}}{7^{4}23^{2}61^{2}}$ & $\frac{2^{8}5^{4}11^{2}17^{2}}{7^{4}23^{2}61^{2}}$ \\
$-680$ & $\frac{11^{3}17^{3}41^{3}}{2^{4}7^{6}23^{2}}$ & $\frac{5^{10}11^{4}}{2^{4}7^{6}23^{2}}$ \\
$-692$ & $\frac{2^{7}17^{6}47^{3}}{5^{4}23^{1}31^{2}53^{2}}$ & $\frac{2^{27}17^{6}31^{2}173^{1}}{5^{10}23^{4}53^{2}}$ \\
$-708$ & $\frac{2^{2}3^{6}5^{6}17^{3}41^{3}53^{3}}{7^{4}13^{2}23^{4}29^{2}37^{2}61^{2}}$ & $\frac{2^{4}11^{4}59^{2}97^{2}157^{2}}{7^{4}23^{4}29^{2}37^{2}61^{2}}$ \\
$-712$ & $\frac{3^{12}47^{1}53^{3}}{2^{7}5^{3}71^{2}}$ & $\frac{19^{4}53^{2}79^{2}173^{2}}{2^{7}5^{6}47^{2}71^{2}}$ \\
$-728$ & $\frac{17^{6}29^{3}53^{3}}{2^{4}5^{2}7^{6}13^{4}71^{2}}$ & $\frac{17^{4}19^{4}59^{4}137^{2}}{2^{4}5^{8}7^{6}13^{2}71^{2}}$ \\
$-740$ & $\frac{2^{4}11^{6}53^{3}}{23^{4}29^{2}37^{1}}$ & $\frac{2^{8}5^{8}11^{4}}{23^{4}29^{2}}$ \\
$-187$ & $\frac{2^{12}3^{6}5^{1}11^{3}}{23^{2}31^{2}37^{2}}$ & $\frac{2^{2}13^{4}71^{2}}{5^{2}23^{2}37^{2}}$ \\
$-760$ & $\frac{3^{3}17^{3}47^{3}}{7^{2}31^{2}71^{2}}$ & $\frac{5^{2}11^{2}13^{2}19^{1}37^{2}}{7^{2}31^{2}71^{2}}$ \\
$-772$ & $\frac{2^{2}3^{6}5^{1}29^{3}41^{3}53^{1}}{7^{4}13^{2}37^{2}61^{2}}$ & $\frac{2^{4}17^{2}113^{2}173^{2}193^{1}}{5^{2}7^{4}53^{2}61^{2}}$ \\
$-195$ & $\frac{2^{12}3^{6}}{13^{2}29^{2}}$ & $\frac{2^{2}5^{6}19^{2}}{13^{2}29^{2}}$ \\
$-788$ & $\frac{2^{5}17^{6}47^{3}59^{3}}{7^{10}13^{4}23^{1}}$ & $\frac{2^{21}13^{2}17^{2}19^{2}59^{4}197^{1}}{5^{6}7^{10}23^{4}}$ \\
$-808$ & $\frac{3^{9}5^{4}41^{3}59^{3}}{2^{2}13^{2}23^{2}47^{2}71^{2}79^{2}}$ & $\frac{11^{6}13^{2}157^{2}197^{2}}{2^{2}5^{2}23^{2}47^{2}71^{2}79^{2}}$ \\
$-203$ & $\frac{2^{27}5^{2}11^{3}}{7^{4}13^{4}37^{2}}$ & $\frac{2^{4}11^{6}79^{2}}{5^{4}7^{4}37^{2}}$ \\
$-820$ & $\frac{2^{2}3^{6}59^{3}}{7^{4}31^{2}37^{2}}$ & $\frac{2^{8}5^{7}}{7^{4}31^{2}}$ \\
$-840$ & $\frac{3^{6}23^{3}53^{3}}{2^{4}7^{2}13^{4}79^{2}}$ & $\frac{5^{4}11^{4}17^{2}19^{2}}{2^{4}7^{2}13^{4}79^{2}}$ \\
$-852$ & $\frac{2^{4}3^{12}5^{2}11^{3}47^{3}59^{3}}{13^{6}23^{6}61^{2}}$ & $\frac{2^{12}11^{6}19^{4}71^{2}193^{2}}{5^{4}13^{4}23^{6}61^{2}}$ \\
$-868$ & $\frac{2^{4}3^{12}5^{2}53^{1}}{7^{4}29^{1}37^{2}}$ & $\frac{2^{8}31^{2}137^{2}197^{2}}{5^{4}7^{4}29^{4}}$ \\
$-872$ & $\frac{11^{6}17^{9}59^{3}}{2^{4}7^{10}29^{4}71^{2}}$ & $\frac{11^{6}19^{2}31^{4}71^{2}173^{2}}{2^{4}5^{6}7^{10}29^{4}}$ \\
$-888$ & $\frac{3^{18}41^{3}47^{3}}{2^{4}5^{2}29^{4}31^{2}79^{2}}$ & $\frac{31^{2}37^{2}59^{4}71^{2}97^{2}}{2^{4}5^{8}29^{4}79^{2}}$ \\
$-227$ & $\frac{2^{33}17^{3}}{13^{6}31^{2}}$ & $\frac{2^{5}17^{4}37^{4}}{5^{6}13^{2}31^{2}}$ \\
$-920$ & $\frac{23^{2}59^{3}}{2^{3}29^{1}47^{2}}$ & $\frac{5^{15}37^{2}}{2^{3}23^{1}29^{4}47^{2}}$ \\
$-932$ & $\frac{2^{6}17^{3}23^{2}41^{3}53^{3}}{5^{2}7^{12}31^{4}}$ & $\frac{2^{12}17^{6}19^{4}53^{2}71^{4}233^{1}}{5^{8}7^{12}23^{4}31^{4}}$ \\
$-235$ & $\frac{2^{6}3^{3}17^{3}}{7^{2}37^{2}47^{1}}$ & $\frac{2^{1}5^{2}11^{2}19^{2}}{7^{2}37^{2}47^{1}}$ \\
$-948$ & $\frac{2^{6}3^{21}59^{3}71^{3}}{5^{2}31^{4}37^{2}47^{2}61^{2}}$ & $\frac{2^{26}19^{4}79^{2}157^{2}}{5^{8}37^{2}47^{2}61^{2}}$ \\
$-952$ & $\frac{3^{12}5^{2}71^{1}}{2^{2}7^{4}23^{1}79^{2}}$ & $\frac{17^{2}113^{2}193^{2}233^{2}}{2^{2}5^{4}7^{4}23^{4}71^{2}}$ \\
\hline
\end{longtable}\end{center}
\end{subsection} \end{section}

\end{document}